\input amstex 
\documentstyle{amsppt}
\magnification 1100
\NoBlackBoxes
\voffset= 1 truecm

\hoffset= -.4 truecm
\hsize= 16.85 truecm
\vsize= 22 truecm

\def\ns{\hskip -.15 truecm}
\def\restobs{Observation 1.1 }
\def\vanlemma{Lemma 1.2 }
\def\rsNp{Theorem 1.3 }
\def\acNp{Theorem 1.3 }
\def\question{Observation 1.4 }
\def\extremal{Remark 1.5 }
\def\KA{Proposition 1.10 }
\def\ava{Theorem 1.21 }
\def\vapn{Proposition 1.22 }
\def\Mukai{Theorem 1.23 }
\def\term{Lemma 1.25 }
\def\ReiderNp{Theorem 1.24 }

\def\BS{[BS] }
\def\Bu{[B] }
\def\CM{[Mu] }
\def\EL{[EL] }
\def\GPone{[GP1] }
\def\GPtwo{[GP2] }
\def\GPthree{[GP3] }
\def\GPfour{[GP4] }
\def\GPfive{[GP5] }
\def\Green{[G] }
\def\GL{[GL] }%
\def\Hbone{[Hb1] }
\def\Hbtwo{[Hb2] }
\def\Ha{[Ht] }
\def\Hoone{[Ho1] }
\def\Hotwo{[Ho2] }
\def\JPW{[JPW] }
\def\Mi{[Mi] }
\def\OP{[OP] }

\def\ns{\hskip -.15 truecm}
\def\seq{(*) }
\def\obsfree{\freeobs}
\def\Spanishlemma{\Splemma}
\def\Splemma{Lemma 2.5 }

\def\freeobs{Observation 2.4 }
\def\genpn{Theorem 1.3 }
\def\tensNp{Theorem 2.6 }
\def\Fanopn{Theorem 2.1 }
\def\GLlemma{Theorem 2.3 }
\def\FanoNp{Theorem 2.7 }
\def\CYpn{Theorem 3.1 }
\def\CYNp{Theorem 3.2 }

\topmatter
\title Some results on rational surfaces and Fano varieties \endtitle
\author Francisco Javier Gallego \\ and \\ B. P. Purnaprajna
\endauthor
\address{Francisco Javier Gallego: Dpto. de \'Algebra,
 Facultad de Matem\'aticas,
 Universidad Complutense de Madrid, 28040 Madrid,
Spain}\endaddress
\email{gallego\@eucmos.sim.ucm.es}\endemail
\address{ B.P.Purnaprajna:
405 Snow Hall,
  Dept. of Mathematics,
  University of Kansas,
  Lawrence, Kansas 66045-2142
}\endaddress
\email{purna\@math.ukans.edu}\endemail
\abstract
The goal of this article is to study the equations and  syzygies of 
embeddings of rational surfaces and certain Fano varieties. Given a rational
surface $X$ and an ample and base-point-free line bundle $L$ on $X$, we
give an optimal numerical criterion for $L$ to satisfy property $N_p$. 
This criterion turns
out to be a characterization of property $N_p$ if $X$ is anticanonical. We
also prove syzygy results for adjunction bundles and a Reider type theorem for
higher syzygies. 

For certain Fano varieties we also prove results on very ampleness and higher
syzygies.

\endabstract
\endtopmatter
\document 
\vskip .3 cm

\headline={\ifodd\pageno\rightheadline \else\leftheadline\fi}
\def\rightheadline{\tenrm\hfil \eightpoint SOME RESULTS ON RATIONAL SURFACES
AND FANO VARIETIES
 \hfil\folio}
\def\leftheadline{\tenrm\folio\hfil \eightpoint F.J. GALLEGO \& B.P. PURNAPRAJNA \hfil}

\heading 0. Introduction \endheading

The goal of this article is to study the equations and the syzygies of 
embeddings of rational surfaces and certain Fano varieties. Previously Butler, 
Homma, Kempf, and the authors had proved results regarding syzygies of
(geometrically) ruled surfaces and surfaces of nonnegative Kodaira dimension.
We will be interested in knowing under what conditions the resolution of the
homogeneous coordinate ring $S/I$ of an embedded variety is ``simple". More
precisely we want to know under what conditions the so-called
property
$N_p$ after M. Green is satisfied. We define this property next:

\proclaim{Definition 0.1} Let $X$ be a projective variety. A very ample line
bundle
$L$ is said to satisfy property $N_0$ if $|L|$ embeds $X$ as a projectively normal
variety. A very ample line bundle
$L$ satisfies property
$N_1$ if
$L$ satisfies property
$N_0$ and the homogeneous ideal $I$ of the image of $X$ embedded by $|L|$ is
generated by quadratic equations. Finally a very ample line bundle $L$ is said
to satisfy property $N_p$, $p \geq 1$,  if it satisfies property $N_1$ and the
matrices in the minimal graded free resolution of $S/I$ have linear entries
from the second  to the
$p$-th step.
\endproclaim

Section 1 is devoted to rational surfaces. 
Given a
rational surface $X$ and an ample and base-point-free line bundle $L$ on $X$, we
observe there is an optimal criterion for $L$ to satisfy  property
$N_p$. This criterion (cf. \rsNp
\ns) depends solely on the intersection number of $L$ with $-K_X$, namely, $L$
satisfies property $N_p$ if $-K_X \cdot L \geq p+3$. This criterion turns out
to be a characterization if $X$ is anticanonical, i.e., a rational surface with
effective anticanonical class. Anticanonical surfaces have been extensively 
studied by several authors, among them B. Harbourne, and \rsNp improves and
generalizes one of his (cf. \Hbtwo \ns). 

\bigskip 

We also 
study the syzygies associated to adjunction bundles and prove a
Reider type theorem for higher syzygies. More precisely,
given any
$A_1,
\cdots, A_n$ ample line bundles on an anticanonical rational surface $X$ of
fixed 
$K_X^2$, we find a sharp bound on $n$ so that 
$K_X+A_1+
\cdots + A_n$ satisfies property $N_p$. One of the easy consequences of this is
Mukai's conjecture, which is in fact optimal for anticanonical surfaces with
$K_X^2=1$. 

\bigskip

Properties such as base-point-freeness and very ampleness are governed 
numerically, as classical results on curves and Reider's theorem on surfaces
show. It is natural to ask whether the same philosophy holds for the 
property $N_p$, which are a natural generalization of base-point-freeness and
very ampleness. For curves Green's theorem (cf. \Green \ns, Theorem 4.e.1),
which says that a line bundle of degree greater than or equal to $2g+p+1$
satisfies property
$N_p$, provides an affirmative answer. For rational surfaces
results in this article, namely, the already mentioned \rsNp and \ReiderNp \ns,
which is a Reider type result for property $N_p$, show us as well that property
$N_p$ only depends  on numerical criteria. In fact, if
$X$ is anticanonical,
\rsNp tells that property $N_p$ for $L$ is exclusively
governed by the intersection number of $L$ with a particular curve lying on $X$, namely,
an anticanonical curve.  This phenomenon can be also observed in
other surfaces such as elliptic ruled surfaces. Indeed, Homma and the authors 
gave in \Hoone \ns, \Hotwo and \GPone a characterization for properties
$N_0$ and
$N_1$ in terms of the intersection number of $L$ with a few curves lying on the elliptic
surface. In this case the relevant curves are a minimal section
and a fiber of the elliptic ruled surface in addition to the anticanonical 
curve. The authors also gave in \GPtwo a criterion for property $N_p$ in
terms of the intersection number of
$L$ with the three above mentioned curves, and  conjectured that those 
intersection numbers should also characterize  property $N_p$ as they did
characterize property
$N_0$ and $N_1$. It is surprising that results of such similar spirit
hold for both rational surfaces and elliptic ruled surfaces. Even though both are surfaces of
Kodaira dimension $-\infty$, they differ in the fact that the Picard
group of a rational surface is discrete, the same does not happen for an
elliptic ruled surface. The differences can be also be seen in the
methods of proof used in \GPone and \GPtwo and in this article, which are indeed
very distinct. 

\bigskip

We also prove a result (cf. Theorem 1.29) connecting property $N_p$, which is
an extrinsic property depending on the embedding of $X$, with the
``termination" of ampleness of $mK_X+L$. In particular we show for what $m$
the line bundle $mK_X+L$ stops being ample for a line bundle $L$ satisfying
property $N_p$. The formula we obtain for such an $m$ depends on $K_X^2$ and
$p$. 

\bigskip 

In Section 1 we also construct several families of examples. These examples
show that all the theorems and propositions proved are sharp, and that the
bounds cannot be reduced.

\bigskip

In Section 2 we study $n$-dimensional Fano
varieties of index greater than or equal to $n-1$. Let $H$ be primitive line
bundle such that $-K_X=iH$. First we prove \Fanopn \ns, which tells exactly
what property $N_p$ is satisfied by $H$. We prove results regarding very
ampleness and on the higher syzygies of multiples of a primitive
$H$ on
$X$ such that $-K_X=iH$. We derive these syzygy
results 
 from the vanishings of
certain Koszul cohomology groups on the Fano variety $X$. We reduce
these vanishings to the vanishings of similar Koszul cohomology groups on
lower dimensional Fano varieties. These lower dimensional are subvarieties of
$X$, but since we need them to be again Fano varieties, we do not obtain them
by taking subsequent hyperplane sections. We take this more indirect approach
because the techniques of Section 1 do not work for these higher dimensional
varieties. The reason is that the
 information available on the resolution of the coordinate ring of the
subsequent hyperplane sections is not good enough. 

\bigskip

Finally in Section 3 we deal with $n$-dimensional Fano varieties of index
$n-3$. We first give a criterion as to when $nH$ satisfies very ampleness and
property $N_0$ when $n \geq 2$. This criterion is actually a characterization
if $n \geq 3$. Then we prove a result on the higher syzygies of multiples of
$H$. The arguments and techniques used are similar to those used in Section 2.
There is though a difference worth noting. When one is considering an
$n$-dimensional Fano variety $X$ of index greater than or equal to $n-1$, one
deduces the vanishing of Koszul cohomology groups on $X$ from the vanishings of
similar cohomology groups on a rational surface. However, if
$X$ has index $n-3$, one deduces the vanishing of Koszul cohomology groups
on $X$ from the vanishings of similar cohomology groups on a Calabi--Yau
threefold, and eventually, on a surface of general type. The reader might
wonder about what the situation is for $n$-dimensional Fano varieties of index
$n-2$. Those surfaces were studied in \GPfive \ns, where the authors obtained
results which are similar to those in Sections 2 and 3 of this article.

\heading 1. Rational surfaces \endheading

In this section we study property $N_p$ for rational surfaces. The first
result is a criterion for property $N_p$ in terms of a very precise numerical
condition, namely the intersection number of the line bundle $L$ under
consideration with the anticanonical class of the surface. 
In order to prove this theorem we will need two preliminary results we mention
now:

\proclaim{\restobs (\GPfour \ns, Observation 2.3)}
 Let $X$ be a regular variety (i.e., a variety such
that 
$H^1(\Cal{O}
_X)=0).$ Let $E$ be a vector bundle on $X$, let $C$ be a
divisor such that  
$L$
$=\Cal{O}_X\left( C\right) $ is  globally generated line bundle and 
$H^1(E\otimes L^{-1})=0.$ If the multiplication map 
$
H^0(E\otimes \Cal{O}_C)\otimes H^0(L\otimes
\Cal{O}_C)\to H^0(E\otimes L\otimes \Cal{O}_C)
$
surjects, then the map 
$
H^0(E)\otimes H^0(L) \to H^0(E\otimes L)
$
also surjects.
\endproclaim

The second result we need is a useful lemma on the vanishing of cohomology of
big and base-point-free line bundles: 

\proclaim {\vanlemma} Let $X$ be a smooth 
surface with $p_g=q=0$ and $B$ a big base-point-free line bundle on
$X$ such that
$-K_X \cdot B >0$. Then 
$h^1(B)=h^2(B)=0$.
\endproclaim 

{\it Proof.} Since $B$ is big and base-point-free, by Bertini, there exists
smooth, irreducible $C \in |B|$. Consider
$$0 \longrightarrow \Cal O_X \longrightarrow B \longrightarrow
B\otimes \Cal O_C \longrightarrow 0 \ .$$
Since $p_g=0$, $h^2(B)=0$. Since $-K_X \cdot B >0$, by adjunction
deg$(B \otimes \Cal O_C) > 2g(C)-2$, and consequently, $h^1(B)=0$, since
$q=0$. $\square$

We are now ready to state and prove the numerical criterion for property
$N_p$. This criterion turns out to be a characterization of property $N_p$ if
the surface is anticanonical. We remark that the case of property $N_0$
was observed
by B. Harbourne (cf. \Hbtwo \ns).

\proclaim{\rsNp} Let $X$ be a rational surface and let $L$ be an
ample line bundle on $X$. If $L$ is base-point-free and 
$-K_X
\cdot L
\geq p+3$, then $L$ satisfies property $N_p$. 
In addition, if $X$ is anticanonical and $L$ is ample, then $L$ satisfies
property $N_p$ if and only if $-K_X \cdot L \geq p+3$.
\endproclaim

\noindent {\it Proof.} First we prove the part of the result stated for
general rational surfaces, and  we start showing that if
$-K_X
\cdot L
\geq 3$, 
$L$ satisfies property
$N_0$. We want to show that 
$$H^0(rL) \otimes H^0(L) @>\alpha>> H^0((r+1)L) \text{ for all }
r  \geq 1.$$
Since $L$ is ample and base-point-free, we can choose a smooth
and irreducible curve $C
\in |L|$.   Then, by \vanlemma and \restobs \ns, the surjectivity of $\alpha$
 follows from the surjectivity of 
$$H^0(rL_C) \otimes H^0(L_C) @>\beta>> H^0((r+1)L_C) \text{ for
all } r  \geq 1.$$ Since $-K_X \cdot L \geq 3$, deg$L_C \geq
2g(C)+1$, hence by Castelnuovo's Theorem 
$\beta$ surjects. This proves that $L$ satisfies property $N_0$.
\medskip

\medskip

Now we prove the result for general $p$. 
We have just proven that if $-K_X \cdot L \geq p+3$, $L$ satisfies
property $N_0$, i.e., $L$ is very ample and embeds 
$X$ as a projectively normal variety. On the other hand,  by
\vanlemma \ns, Kodaira Vanishing Theorem and because $X$ is
regular,
$H^1(rL)=0$ for all
$r \in \bold Z$. Therefore the image $Y$ of $X$ by the embedding
induced by $|L|$ is arithmetically Cohen-Macaulay.  
Let $H$ be a general
hyperplane of
$\bold P^N$ and let $D$ be the corresponding (smooth, irreducible)
divisor on $X$. Since $X$ is regular, $|L_D|$ embeds $D$ in $H$
and the image of this embedding is 
$Y
\cap H$, which is projectively normal. This is the same as saying
that $L_D$ satisfies property $N_0$.  Moreover, since $Y$ is
arithmetically Cohen-Macaulay, 
the minimal resolution of $Y \cap H$ has the same graded Betti as
the minimal resolution of $Y$ (see also \text{\Green \ns,} Theorem 3.b.7). 
Since $-K_X \cdot L \geq
p+3$ we have by adjunction that deg$ L_D \geq 2g+p+1$, then by Green's theorem
(cf. \Green \ns, Theorem 4.a.1) the homogeneous coordinate ring of $Y \cap H$ is
generated by quadrics and the resolution of its homogeneous coordinate ring is
linear until the $p-$th step. Since we did already prove that $L$ satisfied
property
$N_0$, we see that $L$ satisfies property $N_p$.

\medskip

Let now $X$ be anticanonical. On an anticanonical surface a nef line bundle 
$L$ such that $-K_X \cdot L \geq 2$ is base-point-free 
(cf. \Hbone \ns, Theorem III.1). 
Therefore we have just proven one of the implications, namely that if
$-K_X
\cdot L \geq p+3$, then $L$ satisfies property $N_p$. 
We prove now the other implication. Let $L$ be a line bundle on $X$
satisfying property
$N_p$. In particular $L$ is very ample and embeds $X$ as a
projectively normal variety. Moreover, $H^1(rL)=0$ for all $r \in
\bold Z$. Indeed, $H^1(\Cal O_X)=0$ because $X$ is rational and
$H^1(rL)=0$ by 
 Kodaira Vanishing Theorem if $r$ is negative. If $r$ is positive,
$-K_X \cdot rL >0$, because $L$ is ample and $-K_X$ is effective.
Hence $H^1(rL)=0$ by \vanlemma \ns. Therefore the image of the
embedding of $X$ by $|L|$ is arithmetically Cohen-Macaulay. Let
$C$ be a smooth curve in $|L|$. As observed before, $L$ 
satisfies property
$N_p$ if and only if $L_C$ satisfies property $N_p$. Since $-K_X$
is effective and $C$ moves (for $L$ is very ample),
then $-K_X \otimes
\Cal O_C$ is also effective. Hence
$L_C=K_C + N$  with $N$ effective line bundle of
degree $-K_X \cdot L$. Now if $-K_X \cdot L$ were less than or equal
to $p+2$, $L_C$ would not satisfy property $N_p$. 
This follows from a result
of Green and Lazarsfeld (cf. \GL \ns, Theorem 2)
which says in particular
that on a smooth, irreducible, genus
$g$ curve $C$, a  line bundle of the form $K_C + N$ with $N$
effective of degree $p+2$ does not satisfies property $N_p$. 
$\qed$ 

\bigskip

\proclaim {\question } If $X$ is a rational surface which is not
anticanonical and $L$ is a line bundle on $X$ which satisfies $N_p$, in
general, $-K_X \cdot L$ need not be greater than or equal to $p+3$. However,
there are line bundles $L$ such that $-K_X\cdot L=p+3$ and $L$ satisfies $N_p$
and not $N_{p+1}$. Hence $\rsNp$ is sharp for non-anticanonical surfaces.
\endproclaim 

\noindent{\it Proof.} Consider a non-anticanonical rational surface
$X$ and let
$L$ be an ample and base-point-free line bundle so that
$K_X+L$ is also ample, and such that the general curve in $|L|$ is not
hyperelliptic. Then  
$-K_X\cdot L=p+3$ for some $p$. Let $C$ be general and therefore smooth and
irreducible curve in $|L|$. We know by \rsNp that
$L$ satisfies property 
$N_p$. Since $K_X+L$ is ample, by
Kodaira vanishing and duality, $H^1(-K_X-L)=0$. Since $X$ is not anticanonical, 
$-K_X
\otimes \Cal O_C$ is non effective. Then, arguing as in the proof of \rsNp we
conclude, using \GL \ns, Theorem 2, that $L$ satisfies property 
$N_{p+1}$.

Now we show by means of an example the existence of $(X,L)$, where $X$ is a
non-anticanonical rational surface and $L$ is an ample and free line bundle
such that $L$ satisfies $N_p$ but not $N_{p+1}$ and $-K_X\cdot L=p+3$. Let
$\pi: X \longrightarrow \bold F_0$ be the blowing up of $\bold F_0$ at $9$
points. We choose the $9$ points $p_1, \dots p_{9}$ so that 
$X$ is not anticanonical. Let $E_1, \dots, E_{9}$ be the exceptional divisors
and let
$L=\pi^*(2C_0+nf)-E_1-\cdots -E_{9}$ with
$n \geq 4$. Then $L$ is ample, because $L=A+\pi^*((n-3)f)$, where $A$ is as
in Example 1.16. $L$ is also base-point-free. This can be checked using Reider's
Theorem. Indeed, $L$ can be written as $K_X+2A+\pi^*((n-4)f)$.
Let $C$ be a smooth curve in $|L|$. We
want to show that $-K_X \otimes \Cal O_C$ is effective. We look at the long
exact sequence relating the cohomology of $-K_X-L$, $-K_X$ and $-K_X \otimes
\Cal O_C$. By our choice of $p_1, \dots, p_{9}$, $h^0(-K_X)=0$ and by
Riemann--Roch
$h^1(-K_X)=1$. Also by Riemann--Roch, $h^1(-K_X-L)=n-3$. Then $h^0(-K_X \otimes
\Cal O_C) = n-3 \geq 1$. Now $-K_X \cdot L= 2n-5$. Then according to \rsNp
$L$ satisfies property $N_{2n-8}$ and by \GL \ns, Theorem 2, 
$L$ does not satisfy property $N_{2n-7}$. 

\bigskip

\noindent {\bf \extremal \ns.} \acNp is an example of the following
philosophy: On a variety $X$ the failure of a line bundle $L$ to
satisfy property $N_p$ can be traced to the existence of an {\it
extremal} curve $C$ on $X$. The curve
$C$ is extremal with respect to $L$ and property $N_p$ in the following sense:
$L_C$ satisfies property $N_{p}$ but not property $N_{p+1}$. In the proof of
\acNp the existence of such an extremal curve, namely a smooth curve
$C$ in
$|L|$, plays a key role. We would like to point out the existence
of another extremal curve on $X$. If there is a smooth
irreducible curve $C'$ in $|-K_X|$ (for instance, if $X$ is a
Del Pezzo surface), then $C'$ is a smooth elliptic curve. On an
elliptic curve a line bundle satisfies property $N_p$ if and only
if its degree is greater than or equal to $p+3$. This  follows
from Green's theorem (cf. \Green \ns, Theorem 4.a.1) and the self-duality of the
resolution of any elliptic normal curve or alternatively by the theorem of
Green and Lazarsfeld quoted above (cf. \GL \ns, Theorem 2). Now
\acNp says precisely that
$L$ satisfies property $N_p$ but not property $N_{p+1}$ if and only
if $-K_X \cdot L = p+3$. This agrees with the philosophy just stated, since
$-K_X
\cdot L$ is the degree of $L_{C'}$ and $p+3$ is the degree of those
line bundles of $C'$ satisfying property $N_p$ but not property
$N_{p+1}$.

\bigskip

In the remainder of this section we will focus on anticanonical surfaces.
One of our main purposes is  to find uniform and optimal bounds on
$n$ so that a line bundle of the form $K_X + A_1 + \cdots + A_n$ satisfies
property
$N_p$, where
each $A_i$ is an ample line bundle on $X$. The bounds we will obtain will be
for anticanonical surfaces with a fixed value of $K_X^2$. They will 
therefore be finer than a uniform bound valid for any anticanonical surface.
There are two ingredients we need to attack the problem. First we need to
find a sharp, uniform lower bound on
$n$ so that
$K_X + A_1 + \cdots + A_n$ be very ample. This will be dealt with in Proposition
1.6. Second, because of the 
  numerical characterization of property
$N_p$ given by
\rsNp,  we   need to find  a sharp, uniform lower
bound for
$-K_X \cdot A$ for an ample $A$. This is what we do in Proposition 1.9 and \KA
\ns.  The sharpness of Propositions 1.6, 1.9 and 1.10 is shown by the Examples
1.11 to 1.20 presented later on.

\proclaim{Proposition 1.6} Let $X$ be an anticanonical surface, let
$A_1,\dots,A_n$  be  ample line bundles on $X$, and let
$L=K_X+A_1+\cdots+A_n$.  Then $L$ is very ample if
\item{1)}  $K_X^2=9$ and $n \geq 4$.
\item{2)}  $K_X^2=8$ and $n \geq 3$.
\item{3)}  $3 \leq K^2_X \leq 7$ and $n \geq 2$. 
\item{4)}  $K_X^2= 2$ and $n \geq 3$ or $n \geq 2$ unless $n=2, A_1=A_2=-K_X$. 
\item{5)}  $K^2_X=1$ and $n \geq 4$ or $n \geq 2$ unless 
\itemitem{5a)} $n=2$ and $A_1=A_2=-K_X$ or $A_1=-K_X$ and $A_2=-2K_X$ (or
vice versa); or
\itemitem{5b)} $n=3$ and $A_1=A_2=A_3=-K_X$.
\item{6)}  $K_X^2=0$ and $n \geq 3$.
\item{7)}  $K_X^2<0$ and $n \geq 2$. 

\endproclaim

In order to prove the result we need this

\proclaim {Lemma 1.7} Let $X$ be an anticanonical surface and
let $A$ be an ample line bundle on $X$. Then $A^2 \geq 2$
unless $(X,A)=(\bold P^2,\Cal O_{\bold P^2}(1))$ or $K_X^2=1$
and $A=-K_X$. 
\endproclaim

\noindent {\it Proof.} If $-K_X \cdot A \geq 2$, then by \Hbone \ns, Theorem
III.1.a
$A$ is base-point-free, hence if $A^2 =1$, then $(X,A)=(\bold
P^2,\Cal O_{\bold P^2}(1))$. If $-K_X \cdot A =1$ and $A^2=1$,
then
$(K_X+A)
\cdot A=0$ and hence $K_X+A$ is effective by Riemann--Roch.
Since $A$ is ample, $A=-K_X$ and $K_X^2=1$. $\square$

\bigskip

\noindent (1.8) {\it Proof of Proposition 1.6.} First, we note that if $n \geq
4$, then
$L$ is very ample by Reider's Theorem.  Now, except if
$(X,A)=(\bold P^2,\Cal O_{\bold P^2}(1))$ or $A=-K_X$ and
$K_X^2=1$, an  ample line bundle $A$ satisfies
$A^2 \geq 2$ by Lemma 1.7. Therefore if $n \geq 3$, $L$ is
very ample by Reider's Theorem unless all
$A_i=A$ and, either $(X,A)=(\bold P^2,\Cal O_{\bold P^2}(1))$,
or
$A=-K_X$ and $K_X^2=1$.  

Now we prove that under the
hypothesis stated in the proposition, $L$ is very ample if $n
\geq 2$. 
First we show that the result is true if $A_1=A_2=A$ and
$A^2 \geq 3$. By Reider's Theorem the only obstacle to the
very ampleness of $L$ would be the existence of an
irreducible curve $E$ with $A\cdot E=1$ and $E^2=0$. 
Then we prove the following:

\medskip

\noindent (1.8.1) Let $X$ be an anticanonical surface with $K_X^2 \neq 0, 8$ and
let 
$A$ be an ample line bundle on $X$. Then there does not exist  an irreducible
curve $E$ on $X$ with $A \cdot E=1$ and $E^2=0$.  

\medskip

We
show  that, under the above hypothesis, $p_a(E)=0$ or $1$. Since $E$ is
irreducible and
$E^2=0$, $E$ is nef, hence $-K_X \cdot E \geq 0$ and
$(K_X+E)\cdot E \leq 0$. We exclude both possibilities, starting
with
$p_a(E)=1$. In this case
$-K_X
\cdot A =1$, otherwise $A$ would be base-point-free by \Hbone \ns, Theorem III.1
and this would contradict the fact that $A \cdot E=1$. Now
$K_X+E$ is effective by Riemann--Roch. Indeed, $(K_X+E)\cdot
E=0$ and $E$ is a curve, hence $h^0(K_X+E)=1+h^1(K_X+E)
\geq 1$.  On the other hand $(K_X+E) \cdot A=0$, hence
$E=-K_X$ and $K_X^2=0$, which is excluded by hypothesis. 
Now we rule out the second possibility, namely, $p_a(E)=0$. In
this case $E=\bold P^1$ and it follows from Lemma 4.1.10 in
\BS
that $X$ is a $\bold P^1$-bundle, which is also excluded by
hypothesis. 

By Lemma 1.7 the only case left when $A_1=A_2=A$ is $A^2=2$.  Then
$-K_X
\cdot A
\geq 2$ and $A$ is base-point-free by \Hbone \ns, Theorem III.1.a. 
Then $X$ is either $\bold P^1 \times \bold P^1$ or a double
cover of $\bold P^2$ ramified along a smooth quartic, for
which $A=-K_X$ and $K_X^2=2$. Both possibilities are
excluded by hypothesis.

Now we deal with the case $A_1 \neq A_2$. 
If one $A_i$ is such that $A_i^2 \geq 5$, by (1.8.1) and Reider's Theorem
$K_X+A_1+A_2$ is very ample. If, say, $A_1^2=4$ and $A_2^2 \geq 2$, we prove
the following:

\noindent (1.8.2) Let $X$ be an anticanonical surface and let $A_1$ and
$A_2$ be ample line bundles on $X$ such that $A_1^2 \geq 3$. Then $A_1
\cdot A_2
\geq 2$. 

Then we are  done by (1.8.1), (1.8.2) and Reider's Theorem. Therefore we prove
(1.8.2). Assume the contrary. Now for any ample bundle $A$ on $X$, and in
particular for $A_2$, we have that $h^2(A) = h^0(K_X-A)=0$, because any ample 
bundle intersects $K_X-A$ strictly negatively. Then by Riemann--Roch $A_2$ is
effective. On the other hand $(K_X+2A_1)\cdot A_2 \leq 1$. 
We have seen before that $K_X+2A_1$
is very ample, therefore $(K_X+2A_1)\cdot A_2 = 1$, $-K_X \cdot A_2=1$  and each
member of
$|A_2|$ is a smooth $\bold P^1$. But if $-K_X \cdot A_2=1$ then 
$(K_X+A_2)\cdot A_2
\geq 0$ and this contradicts the fact that the sectional genus of $A_2$ is $0$.

Now if, say, $A_1^2=4$ and $A_2^2=1$, again we know by (1.8.2) that $A_1 \cdot
A_2 \geq 2$, and if $A_1 \cdot A_2 \geq 3$ we are again done by (1.8.1) and
Reider's Theorem. Thus we consider the case $A_1^2=4$, $A_2^2=1$ and $A_1
\cdot A_2 =2$. Since $n =2$ and $A_2^2=1$, by Lemma 1.7 $A_2=-K_X$ and
$K_X^2=1$. 
Then $A_1 \cdot (2K_X + A_1)=0$. On the other hand
$h^2(2K_X+A_1)=h^0(-K_X-A_1)=0$, because $A_1 \cdot (-K_X-A_1)=-2$. Then it
follows by Riemann-Roch that $2K_X + A_1$ is effective, so $A_1=-2K_X$.
However the possibility $K_X^2=1$, $A_2=-K_X$ and $A_1=-2K_X$ is excluded by
hypothesis. 

Now the only remaining cases are $A_1^2=3$ and $A_2^2=1,2,3$, $A_1^2=2$ and
$A_2^2=2,1$ and $A_1^2=A_2^2=1$. We analyze them  case by case. If 
$A_1^2=A_2^2=3$, then by (1.8.2) $A_1 \cdot A_2 \geq 2$, then by (1.8.1) and
Reider's Theorem we are done. If $A_1^2=3$ and $A_2^2=2$, then $-K_X \cdot A_2
\geq 2$ and as already argued $A_2=-K_X$ and $K_X^2=2$. If $A_1^2=3$ and
$A_2^2 =1$, then by Lemma 1.7, $A_2=-K_X$ and $K_X^2=1$. But then in neither of
the cases can it happen that
$A_1^2=3$ and
$-K_X
\cdot A_1=2$. If
$A_1^2=2$ and $A_2^2=2$ or $1$, then either $A_1=A_2=-K_X$ and $K_X^2=2$, which
is excluded because we are assuming $A_1 \neq A_2$, or $A_1=-K_X$ and $K_X^2=2$
and $A_2=-K_X$ and $K_X^2=1$, which is absurd. Finally if $A_1^2=A_2^2=1$, then
by Lemma 1.7, $A_1=A_2=-K_X$ and $K_X^2=1$, but we are assuming $A_1 \neq A_2$.
 $\square$

\bigskip

We now proceed to study lower bounds for $-K_X \cdot A$. One lower bound is of
course $1$. This bound is sharp if $K_X^2 \leq 1$, as proven by Examples 1.15,
1.18, 1.19 and 1.20. However, if $K_X^2 \geq 2$ a better bound holds:  First, if
$K_X^2=9$, then $X = \bold P^2$ and the bound is $3$. If $K_X^2=8$, $X$
is a Hirzebruch surface and we summarize the bounds on $-K_X \cdot A$ in the
next proposition. Finally if $1 \leq K_X^2 \leq 7$, $-K_X \cdot A \geq
K_X^2$. This bound is sharp as shown in Examples 1.13 and 1.14. We end this
study by classifying the boundary and near-boundary cases when $1 \leq K_X^2
\leq 7$  and by obtaining a better bound when $(X,A)$ does not fall in one of
these boundary and near-boundary cases. As explained before, we need these
technical results to prove \text{\Mukai \ns.} 

\proclaim{Proposition 1.9} Let $X$ be an anticanonical surface and let $A$ be
an ample line bundle on $X$.
If $K_X^2=8$, i.e., if $X$ is a (geometrically) ruled rational
surface, then $-K_X \cdot A \geq 4$. More precisely, if $X =
\bold F_e$, then
$-K_X
\cdot A
\geq e+4$.
\endproclaim

\noindent {\it Proof.} If $X = \bold F_e$,
$-K_X$ is linearly equivalent to $2C_0 + (e+2)f$, where $C_0$ is the minimal
section of
$\bold F_e$ and
$f$ is a fiber. In any case there is a divisor in the anticanonical class
having $e+4$ irreducible components counted with multiplicity. Then the
intersection number of any ample line bundle $A$ with this divisor is greater
than or equal to $e+4$. $\square$

\proclaim{\KA } Let $X$ be an anticanonical surface and
let $A$ be an ample line bundle such that 
 $1 \leq K_X^2 \leq 7$.

Then $-K_X \cdot A \geq K_X^2 +3$ unless one of the following happens:

\item{a)} $A=-K_X$, in which case $-K_X \cdot A = K_X^2$;

\item{b)} $K_X^2=1$ and $A=-2K_X$, in which case $-K_X \cdot A = K_X^2 +1$;

\item{c)} $K_X^2=1$ and $A=-3K_X$, in which case $-K_X \cdot A = K_X^2 +2$;

\item{d)} 
$K_X^2=2$ and
$A=-2K_X$, in which case
$-K_X\cdot A=K_X^2+2$;

\item{e)} 
$K_X+A$ is a base-point-free line bundle, $A$ is very ample 
and $(X,A)$ is a conic fibration under $|K_X+A|$, in which case
$-K_X\cdot A \geq  K_X^2+2$. 
\endproclaim 

\noindent {\it Proof.} We assume throughout the proof that $A \neq -K_X$ and we
want to see that, except for the other exceptions listed in the statement, 
$-K_X\cdot(K_X+A)
\geq 3$. We divide the proof in two cases:

\noindent {\it Case 1:} $K_X^2 \geq 2$. Assume first that $A^2 \geq 5$. Then by
Hodge Index Theorem
$-K_X
\cdot A \geq 4$; in particular, $A$ is very ample. We apply
Reider's Theorem to $K_X+A$ to see it is base-point-free. The only obstruction
to
$K_X+A$ being base-point-free is the existence of a reduced and irreducible
curve $E$ such that $A \cdot E=1$ and $E^2=0$. Since $A$ is very
ample, $E = \bold P^1$. Then it follows from \BS \ns, Lemma 4.1.10
that $X$ would be a $\bold P^1$-bundle, which is excluded by hypothesis.
We apply now \Hbone \ns, Lemma II.6.a. 
If $K_X+A$ is composed with a
pencil, since $K_X+A$ is base-point-free and $A$ is ample,
$(K_X+A)^2=0$, $h^1(K_X+A)=0$ and $-K_X \cdot (K_X+A) = 2r$, for
some $r \geq 1$. If $r \geq 2$, we are done. If $r=1$, then
$(K_X+A)^2=0$ and $-K_X\cdot (K_X+A)=2$. Hence $(K_X+A)\cdot A=2$,
therefore the sectional genus of $A$ is $2$. 

Then by \BS \ns,
Theorem 10.2.7.2 $(X,A)$ is a conic fibration over $\bold P^1$
under $|K_X+A|$. Note that in this case $-K_X\cdot(K_X+A)$ is an even number
greater than or equal to $2$.

If $K_X+A$ is not compose with a pencil then $(K_X+A)^2 >0$ by \Hbone \ns,
Lemma II.6.b. We want to see under what conditions $-K_X\cdot(K_X+A) \geq 3$.
Since $K_X^2 \geq 2$ it follows from Hodge Index Theorem that
$-K_X\cdot(K_X+A) \geq 2$. Assume then that $-K_X\cdot(K_X+A) =
2$. Then $(K_X+A)^2 \geq 2$ and even. 
Now if $(K_X+A)^2 \geq 4$ or $K_X^2 \geq 3$, we are done by Hodge Index
Theorem. Then the only case left is $(K_X+A)^2=K_X^2=2$. 

On the other hand $K_X+A$ is big and nef, therefore
by Kawamata-Viehweg $h^1(2K_X+A)=h^2(2K_X+A)=0$ and by
Riemann--Roch $2K_X+A$ is effective.
 Since $A$ is ample and $A \cdot
(2K_X+A) = 0$,
this implies that $A=-2K_X$. 

We deal now with the case $A^2 \leq 4$. By Lemma 1.7, $A^2 \geq 2$,
and by Hodge Index Theorem $-K_X \cdot A \geq 2$, hence by \Hbone \ns, Theorem
III.1.a
$A$ is base point free. If $A^2=2$, then either $X= \bold P^1 \times \bold
P^1$ or $K_X^2=2$ and $A=-K_X$. The two  possibilities are  excluded either by
hypothesis or by assumption. 
Let now $A^2=3$ or
$4$. 

If $A^2 = 3$, by Hodge Index
Theorem
$-K_X
\cdot A
\geq 3$, hence by \rsNp \ns, $A$ is very ample. Then $X$ is a
 Del Pezzo cubic surface in $\bold P^3$ and $A = -K_X$ or a
rational normal scroll in $\bold P^4$. The   two possibilities are excluded by
hypothesis or by assumption.

Finally, if $A^2=4$, again by Hodge Index Theorem $-K_X
\cdot A \geq 3$ and $A$ is very ample. However the only linearly
normal smooth surfaces of degree $4$ are K3 surfaces in $\bold
P^3$, the Del Pezzo surface in $\bold P^4$ and the Veronese
Surface and the rational normal scrolls in $\bold P^5$. 
\medskip

\noindent {\it Case 2:} $K_X^2 = 1$. Assume first that $A^2 \geq 5$. Then by
Hodge Index Theorem
$-K_X
\cdot A \geq 3$; in particular, $A$ is very ample. The same argument used in
Case 1 proves that $K_X+A$ is base-point-free.  We apply again \Hbone \ns, Lemma
II.6.a. 
Then if $K_X+A$ is composed with a
pencil, we get as before that either $-K_X \cdot (K_X+A) \geq 4$ or, $-K_X
\cdot (K_X+A) = 2$ and by
\BS
\ns, Theorem 10.2.7.2 $(X,A)$ is a conic fibration over $\bold P^1$
under $|K_X+A|$. 

If $K_X+A$ is not compose with a pencil then $(K_X+A)^2 >0$ by \Hbone \ns,
Lemma II.6.b. We want to see under what conditions $-K_X\cdot(K_X+A) < 3$.
If $(K_X+A)^2 \geq 5$, then $-K_X\cdot(K_X+A) \geq 3$ by Hodge Index Theorem.
Then we have to study the cases $1 \leq (K_X+A)^2 \leq 4$. Then by Hodge Index
Theorem $-K_X\cdot(K_X+A) \geq 1$. 

Assume $(K_X+A)^2=1$. Then $-K_X\cdot(K_X+A)$ is odd, so we only have to worry
about $-K_X\cdot(K_X+A)=1$. Recall that $K_X+A$ is free and big, therefore by
Kawamata--Viehweg Theorem, $h^1(2K_X+A)=h^2(2K_X+A)=0$. Since $(2K_X+A)\cdot
(K_X+A)=0$, by Riemann--Roch it follows that $h^0(2K_X+A)=1$. On
the  other hand, since $K_X^2=(K_X+A)^2=1$, then $A\cdot (2K_X+A)=0$, hence
$A=-2K_X$. This is excluded because at this point we are assuming $A^2 \geq 5$. 

Assume $(K_X+A)^2=2$. Then $-K_X\cdot(K_X+A)$ is even, so we only have to worry
about $-K_X\cdot(K_X+A)=2$. Note that $A\cdot (K_X+A)=4$. Then by Riemann-Roch
$|K_X+A|$ maps $X$ as double cover of $\bold P^2$. Then either $K_X^2=8$, which
is excluded by hypothesis or $K_X^2=2$, which is excluded by assumption. 

Assume $(K_X+A)^2=3$. Then $-K_X\cdot(K_X+A)$ is odd, so we only have to worry
about $-K_X\cdot(K_X+A)=1$. This cannot happen by Hodge Index Theorem. 

Assume $(K_X+A)^2=4$. Then $-K_X\cdot(K_X+A)$ is even, so we only have to worry
about $-K_X\cdot(K_X+A)=2$. First we see that $(2K_X+A)$ is effective. Indeed,
$(2K_X+A)\cdot(K_X+A)=2$, then by Kawamata--Viehweg Theorem and Riemann--Roch,
$h^0(2K_X+A)=2$. Then $h^2(3K_X+A)=h^0(-2K_X-A)=0$. We also have that $(3K_X+A)
\cdot (2K_X+A) =0$. Then by Riemann--Roch $3K_X+A$ is effective. On the other
hand $A \cdot (3K_X+A) \leq 0$, hence $A \cdot (3K_X+A)=0$ and $A = -3K_X$.

We deal now with the case $A^2 \leq 4$. By Lemma 1.7, if $A^2 = 1$, then
$A=-K_X$, which is excluded by assumption. If $A^2 \geq 2$,  by Hodge Index
Theorem
$-K_X
\cdot A
\geq 2$, hence by
\Hbone
\ns, Theorem III.1.a
$A$ is base point free. If $A^2=2$, then either $X= \bold P^1 \times \bold
P^1$, which is excluded by hypothesis or $K_X^2=2$ and $A=-K_X$, which is
excluded by assumption. 

If $A^2 = 3$, by Hodge Index
Theorem and because $-K_X \cdot A$ must be odd
$-K_X
\cdot A
\geq 3$, hence by \rsNp \ns, $A$ is very ample. Then $X$
is a
 Del Pezzo cubic surface in $\bold P^3$ and $A = -K_X$ or a
rational normal scroll in $\bold P^4$. Both possibilities are
excluded either by assumption or by hypothesis.

Finally let $A^2=4$. Again by Hodge Index Theorem $-K_X
\cdot A \geq 2$ and $A$ is base-point-free. If $-K_X \cdot A \geq 3$, $A$
is actually very ample by \rsNp \ns. However the only linearly normal smooth
surfaces of degree
$4$ are K3 surfaces in $\bold P^3$, the Del Pezzo surface in $\bold P^4$ and the
Veronese Surface and the rational normal scrolls in $\bold P^5$. 
Then the only possibility left to study is $A^2=4$ and $-K_X
\cdot A = 2$.
Then $A  \cdot (2K_X + A)=0$. On the other hand
$h^2(2K_X+A )=h^0(-K_X-A )=0$, because $A  \cdot (-K_X-A )=-2$. Then it
follows by Riemann-Roch that $2K_X \cdot A_1$ is effective, so $A =-2K_X$. 
$\square$

\bigskip

Before stating and proving Theorems 1.23 and 1.24, which are consequences of
Propositions 1.6, 1.9 and 1.10 it is important to know that these propositions
are sharp. We do this by means of the following series of examples. Then these
examples will also imply the optimality of our theorems on rational surfaces.

\proclaim{Example 1.11} $K_X^2=9$, $A$ ample, $-K_X\cdot A = 3$, $K_X+3A$ not
ample.
\endproclaim

Of course the only rational surface with $K_X^2=9$ is $X=\bold P^2$, $A=\Cal
O_{\bold P^2}(1)$ attains the bound $-K_X\cdot A \geq 3$ and $K_X+3A=\Cal
O_{\bold P^2}$, hence not very ample. 

\proclaim{Example 1.12} $K_X^2=8$, $A$ ample, $-K_X \cdot A = e+4$, $K_X+2A$
not  ample.  
\endproclaim

The rational surfaces with $K_X^2=8$ are the Hirzebruch surfaces $\bold F_e$.
With these examples we show that the bounds computed in Proposition 1.9 for
$-K_X
\cdot A$ in terms of
$e$ are sharp. If $X=\bold F_e$, let $C_0$ be the minimal section and $f$ a
fiber. The divisor
$A=C_0+(e+1)f$ is ample and $-K_X \cdot A$ achieves the bound $e+4$. Moreover
$K_X+2A=ef$ is free but not ample, hence $A$ provides an example of $K_X+2A$
not being very ample.

\proclaim {Example 1.13} $3 \leq K_X^2 \leq
7$, $A$ ample, $-K_X\cdot A = K_X^2$, $K_X+A$ not ample.
\endproclaim 

We obtain $X$ with $K_X^2=9-i$
by blowing up 
$\bold P^2$ at $2 \leq i \leq 6$ points, which we choose in sufficiently
general position  (not
$3$ of them on a line, not $6$ of them on a conic). This is a
Del Pezzo surface  and $-K_X$ is very ample 
(cf. \Ha \ns, Theorem V.4.6). 
Then 
$(X,-K_X)$ achieves the bound $A^2 \geq K_X^2$ and also provides an example
where
$K_X+A$ is not very ample.

\proclaim {Example 1.14} $K_X^2=2$, $A$ ample, $-K_X\cdot A = K_X^2$ and
$K_X+2A$ not very ample.
\endproclaim

Let  $p: X \longrightarrow \bold P^2$ be the double cover of $\bold P^2$ ramified
along a smooth quartic. 
Such an $X$ is rational, 
$K_X=p^*(\Cal O_{\bold P^2}(-1))$ and
$-K_X=p^*(\Cal O_{\bold P^2}(1))$. Therefore $K_X^2=2$ and, since  $p$ is
finite, $-K_X$ is ample and base-point-free. 
In fact, $p$ is induced by the morphism of the complete
anticanonical linear series. We remark that, actually, a rational surface $X$ 
with
$K_X^2=2$ and with anticanonical divisor ample 
is a double
cover of $\bold P^2$ as the one described above. This follows from
Riemann--Roch.  Thus the anticanonical divisor on such a surface $X$ attains
the bound $A^2 \geq K_X^2$. 

Such a surface $X$ can   also be found embedded in $\bold P^6$. Indeed,
 $-2K_X$ is very ample, and embeds $X$ as a degree $8$, sectional genus $6$,
smooth surface in $\bold P^6$   (see \BS \ns, Example 10.2.4). 

In addition $A=-K_X$ provides an example of $K_X+2A$ not being very ample.

\proclaim {Example 1.15} $K_X^2=1$, $A$ ample, $-K_X \cdot A=1$, $K_X+3A$ not
very ample.
\endproclaim 

An example of a surface $X$ with $K_X^2=1$ and $-K_X$ ample can be found in
\BS
\ns, Example 10.4.3. 
Indeed there exists a
smooth rational surface with $K_X^2=1$ and $-K_X$ ample. This surface is embedded
by $|-3K_X|$ as a degree $9$,
sectional genus $4$ smooth surface in $\bold P^6$.
Then $(X,-K_X)$ achieves the bound $A^2 \geq K_X^2$. In addition $A=-K_X$
provides an example of $K_X+3A$ not being very ample. Indeed, $-2K_X$ is not
very ample as its complete linear series induces a double cover of quadric
cone in $\bold P^3$.

\proclaim {Example 1.16} $-1 \leq K_X^2 \leq 8$, $A$ ample, $-K_X \cdot A =
K_X^2+2$,
$K_X+A$ not  ample. 
\endproclaim

The family of examples we consider now are conic bundles. Let us fix
$n=K_X^2$, $-1 \leq n \leq 8$.  The examples are constructed from the
pair
$(Y=\bold F_e, 2C_0+mf)$, where
$m=e+3$, $0 \leq e \leq 2$, 
$C_0$ is the minimal section of $Y$ and $f$ is a fiber of $Y$. Let
$C$ be a smooth anticanonical divisor on $Y$. Let $l=8-n$. We choose
$\Sigma=\{p_1,\dots,p_l\}$, $l$ distinct points on $C$ lying on different
fibers of $Y$. Let
$\pi: X \longrightarrow Y$ be the blowing-up of $Y$ along $\Sigma$ and let
$E_1,\dots,E_l$ be the exceptional divisors. Let
$A=\pi^*(2C_0+mf)-E_1-\cdots-E_l$. By the choice of $m$, $K_X+A=\pi^*f$.
Hence $K_X+A$ is not ample and
$-K_X
\cdot A= K_X^2+2$. We see now that $A$ is ample using
Nakai--Moishezon's criterion. Firstly, as $n \geq -1$, $A^2 \geq 3$. Secondly,
it is clear that
$A \cdot E_i$=1, for all $1 \leq i \leq l$. Finally we will check that the
intersection of $A$ with any irreducible non-exceptional curve $T$ on $X$ is
strictly positive. Let $D=\pi(T)$ and let $D \sim aC_0+bf$. Let $m_1, \dots,
m_l$ be the multiplicities of $D$ at $p_1,\dots, p_l$. Then $A \cdot T =
(2C_0+mf)\cdot D -m_1-\cdots -m_l$. Hence we want $m_1+\cdots+m_l <
(2C_0+mf)\cdot D = 3a + 2b -ae$. We first consider the case when $C$ and $D$
intersect properly. Since
$p_1,\dots, p_l \in C$, $m_1+\cdots+m_l \leq C \cdot D$. Therefore 
$m_1+\cdots+m_l \leq 2a+2b -ae$. Then if $a > 0$ we are done. If $a=0$,
then $D=f$. Since we have chosen $p_1,\dots,p_l$ in different fibers in
this case $m_1+\cdots + m_l=1 <2=3a+2b-ae$, and we are also done. Now
consider the case when $C$ and $D$ do not intersect properly. Since both
$C$ and
$D$ are irreducible, $C=D$. In this
case $m_1+\cdots+m_l=l$, for $C$ is smooth, and $a=2$ and $b=e+2$. Then
$3a+2b-ae=10$ and we are done if $l <10$. The latter happens because 
$n
\geq -1$.

\proclaim {Example 1.17} $-2 \leq K^2_X \leq 8$, $A$ ample, $-K_X \cdot A =
K_X^2 + 3$ and
$K_X+A$ not ample. 
\endproclaim

Let $Y=\bold F_1$, let $C_0$ the minimal section, let $f$ be a fiber and let
$C$ be a smooth irreducible anticanonical curve. Let
$\pi:X \longrightarrow Y$ the blowing up of $Y$ at $\Sigma=\{p_1, \dots,
p_l\}$, where $p_1, \dots,p_l$ are distinct points of $Y$ on $C$ away from
$C_0$ and
$0 \leq l
\leq 10$. Let
$E_1,\dots,E_l$ be the exceptional divisors lying over $p_1,\dots,p_l$
respectively. Let $A=\pi^*(3C_0+4f)-E_1-E_2-\cdots-E_l$. We claim that
$A$ is ample, that $-K_X\cdot (K_X+A)=3$ and that $K_X+A$ is not ample. The
latter is clear, since $K_X+A=\pi^*(C_0+f)$. From this it also follows that
$-K_X \cdot (K_X+A)=3$. We check finally the ampleness of A using
Nakai-Moishezon's criterion. On the one hand $A^2 = 15 -l \geq 5$, since
$l \leq 10$. Clearly
$A
\cdot E_i
= 1$. Let
$T$ be now a nonexceptional irreducible curve on $X$, and let $D=\pi(T)$. Let
$D \sim aC_0+bf$ and let $m_1,\dots,m_l$ be the multiplicities of $D$ at
$p_1,\dots,p_l$. Then
$A\cdot T = (3C_0+4f)\cdot D  -m_1  - \cdots - m_l= a+3b -m_1 - \cdots -
m_l$, so we want $m_1+m_2+\cdots+m_l < a+3b$.  We distinguish several
cases: first we consider the case when $C$ and
$D$ intersect properly. Then $m_1 + \cdots +m_l \leq C\cdot D = a+2b$
since every
$p_i$ lies on $C$. We consider two subcases: $b > 0$ or $b=0$. If the
former, $a + 2b < a + 3b$ and we are done. If the latter $T=C_0$ and
$m_1=\cdots=m_l=0$ by our choice of $p_1,\dots,p_l$. Since $0=m_1+\cdots+m_l <
1$  we are also done. The only case left is when $C$ and $D$ do not
intersect properly. Since both
$C$ and $D$ are irreducible, $C=D$. Since $C$ is smooth, $m_1+\cdots+m_l=l$.
On the other hand
$a+3b = 11$ in this case, and, since $l \leq 10$, we are done.

\proclaim {Example 1.18} $K_X^2=0$, $A$ ample, $-K_X \cdot A = 1$, $K_X+2A$ not
very ample. 
\endproclaim

We find a surface with $K_X^2=0$  with an ample line bundle $A$ such that $-K_X
\cdot A=1$ and $K_X+2A$ is not very ample. Let $Y=\bold P^2$. We consider a set
$\Sigma
\subset Y$ of $9$ distinct points being the complete intersection of $2$ cubics, and
neither $3$ of them on a line nor $6$ of them on a conic. Let $X$ be the blowing up
of $Y$ along $\Sigma$. Then $-K_X$ is base-point-free, $h^0(-K_X)=2$, and all $C \in
|-K_X|$ are irreducible curves. Therefore $|-K_X|$ turns $Y$ into an elliptic
fibration $\varphi: X \longrightarrow \bold P^1$ with irreducible fibers and at least
$9$ sections, namely, the
$9$ exceptional divisors. Let $E$ be one of the exceptional divisors and let $F$ be
a fiber of  $\varphi$. Then $A=E + 2F$ is ample. Indeed, we use Nakai-Moishezon's
criterion. The self-intersection $(E+2F)^2 = E^2 + 4 E \cdot F = 3$. If $F'$ is a
fiber of $\varphi$, $(E+2F)\cdot F'=1$. Let $C$ be an irreducible curve on
$X$, and not a fiber. Then either $C=E$, in which case $(E+2F)\cdot C=1$ or
$(E+2F) \cdot C = E\cdot C + 2 F \cdot C \geq 2$. Now we see that $-K_X
\cdot A=1$  and that
$K_X+2A$ is not very ample. It is
clear that
$-K_X
\cdot A=1$. We see  now that
$K_X+2A$ is not very ample. Indeed, let $C \in |-K_X|$. Then $C$ has
arithmetic genus
$1$; however
$K_X+2A=-3K_X+2E$, hence
$(K_X+2A)\cdot C=2$. Since the arithmetic genus of $C$ is $1$, $K_X+2A$ is not very
ample. Note that $A'=E+nF$ for some $n \geq 2$ satisfies as well that $-K_X\cdot
A'=1$ and $K_X+2A'$ is not very ample.

\proclaim {Example 1.19} $K_X^2<0$ odd, $A$ ample, $-K_X \cdot A=1$ and $K_X+A$
not  ample.
\endproclaim

\noindent Given $n$ odd number strictly smaller than $0$, we find a surface
$X$ with
$K_X^2=n$ and an ample line bundle $A$ on $X$ such that $-K_X \cdot A =1$.
Let $Y=\bold F_0$, let $l=8-n$. Let $C$ be a smooth anticanonical curve on
$Y$ and $f_1$ and $f_2$ be two lines each belonging to one ruling of $\bold
F_0$. Let
$\Sigma=\{p_1,\dots,p_l\}$ be $l$ distinct points on $C$ chosen so that not two
of them belong to a line of $\bold F_0$. Let
$\pi: X
\longrightarrow Y$ be the blowing up of $Y$ along $\Sigma$ and let $E_i$ be the
exceptional divisor over $p_i$. Let $k=\frac {l-3} {2}$ and let
$A=\pi^*(\Cal O_{\bold F_0}(2f_1+kf_2)-E_1 -\cdots -E_l)$. We see that $-K_X
\cdot A = (2f_1+2f_2)\cdot (2f_1+kf_2) + E_1^2+\cdots+ E_l^2 =1$. Moreover
$K_X+A=\pi^*((k-2)f_2)$, and therefore it is not ample. Now we will see that $A$
is ample using Nakai--Moishezon's criterion. First
$A^2=l-6
\geq 3$. Now we see the intersection of $A$ with irreducible curves on $X$.
The intersection of $A$ with each $E_i$ is $1$. Let $T$ be an irreducible
curve which is not an exceptional divisor and let $D=\pi(T)$. Let $D \sim
af_1+bf_2$. Let $m_1, \dots, m_l$ be the multiplicities of $D$ at
$p_1,\dots,p_l$. Then $A\cdot T = (2f_1+kf_2)\cdot D -m_1-\cdots -m_l$. Hence
we want to see that $m_1+\cdots +m_l <ka+2b$. Since
$p_1,\dots, p_l \in C$, $m_1+\cdots+m_l \leq C \cdot D$ if $C$ and $D$ intersect
properly. Then
$m_1+\cdots+m_l \leq 2a+2b$. Since $n \leq -1$, $k \geq 3$, hence if $a
\geq 1$, $2a+2b < ka+2b$ and we are done. If $a=0$, then $D=f_2$ and by our
choice of $p_1,\dots,p_l$, $m_1+\cdots + m_l =1 <2= ka+2b$ and we are also
done. Now suppose that
$C$ and
$D$ do not intersect properly. Since both
$C$ and
$D$ are irreducible,  $C=D$ and $a=b=2$. In this
case
$m_1+\cdots+m_l=l$, for
$C$ is smooth.  
Since
$k=\frac {l-3} {2}$, $l < 2k+4$. 

\proclaim{Example 1.20} $K_X^2<0$ even, $A$ ample, $-K_X \cdot A=1$ and $K_X+A$
not  ample.
\endproclaim

\noindent Given $n$ even number strictly smaller than $0$, we find a surface
$X$ with
$K_X^2=n$ and an ample line bundle $A$ on $X$ such that $-K_X \cdot A =1$.
Let $Y=\bold F_0$ and let $l=8-n$. Let $C$ be a smooth anticanonical curve on
$Y$ and $f_1$ and $f_2$ be two lines each belonging to one ruling of $\bold
F_0$. Let
$\Sigma=\{p_1,\dots,p_l\}$ be $l$ distinct points on $C$ chosen so that
not two of them belong to a line of $\bold F_0$. Let
$\pi: X
\longrightarrow Y$ be the blowing up of $Y$ along $\Sigma$ and let $E_i$ be the
exceptional divisor over $p_i$. Let $k=\frac {l-4} {2}$ and let
$A=\pi^*(\Cal O_{\bold F_0}(3f_1+kf_2)-2E_1 -E_2 -\cdots -E_l$. Then $-K_X
\cdot A = (2f_1+2f_2)\cdot (3f_1+kf_2) + 2E_1^2 +E_2^2+\cdots+ E_l^2 =1$. In
addition $K_X + A = \pi^*(f_1+(k-2)f_2)-E_1$. Therefore $K_X+A$ is not
ample, since its intersection with the strict transform of one of the lines
passing through $p_1$ is $0$. We see now that
$A$ is ample  using Nakai--Moishezon's criterion. First
$A^2 = 2l-15
\geq 5$. Now we see the intersection of $A$ with the irreducible curves on
$X$. The intersection of
$A$ with $E_1$ is $2$ and with the  other exceptional divisors is $1$.  Let
$T$ be an  irreducible curve which is not an exceptional divisor and let
$D=\pi(T)$. Let $D
\sim af_1+bf_2$. Let $m_1, \dots, m_l$ the multiplicities of $D$ at $p_1,\dots,p_l$.
Then $A\cdot T = (3f_1+kf_2)\cdot D
-2m_1-m_2-\cdots -m_l$. Hence we want to see that $2m_1+m_2+\cdots +m_l <ka+3b$. 
First we suppose that $C$ and $D$ intersect properly. Since
$p_1,\dots, p_l \in C$, then $m_1+\cdots+m_l \leq C \cdot D=2a+2b$. We
distinguish now two cases. First, if $m_1 \leq a$, then $2m_1+m_2+\cdots +m_l
\leq 3a+2b$. Since $n \leq -2$, $k \geq 3$. Then if $b \geq 1$,  $3a + 2b <
ka + 3b$ and we are done. If $b=0$, then $D \sim f_1$ and $a=1$. By the
choice of
$p_1,\dots,p_l$,
$2m_1+m_2+\cdots+m_l
\leq 2 < 3 \leq ka+3b$. 
Second, if $m_1 > a$, then $D$ passes through $p_1$ and $D \sim  f_2$.
Then,    because of
the choice of
$p_1,\dots,p_l$, $2m_1+m_2+\cdots+m_l \leq 2 < 3=ka+3b$. 
Now we suppose that $C$ and $D$ do not intersect properly.
Since
$C$ and $D$ are irreducible, then $C=D$, and $a=b=2$. In
this case $2m_1+m_2+\cdots+m_l=l+1$, for $C$ is smooth. 
Then we are done if   $l+1 < 2k+6=ka+3b$. This occurs because $k=\frac{l-4}
{2}$.

\bigskip

The lower bounds for $-K_X \cdot A$ obtained in Proposition 1.9 and \KA
combined with \rsNp yield several results. The following generalizes a well
 known  result, namely, the equivalence of the notion of ampleness and
very ampleness for $\bold P^2$ and for Hirzebruch surfaces (cf. \Ha \ns,
Corollary V.2.18).  
The result we get is actually stronger than this classical result since we
obtain the equivalence of ampleness and certain $N_p$ property for certain
anticanonical surfaces:

\proclaim{\ava } Let $X$ be a  rational  surface 
such
that    
$d=K_X^2 \geq 3$.
Let $A$ be
a line bundle on $X$. The following are equivalent: 

\item{1)} $A$ is ample;
\item{2)} $A$ is very ample;
\item{3)} $A$ satisfies property $N_0$. 

More precisely:

\item{--} If $3 \leq K_X^2 \leq 7$, then $A$ satisfies property $N_{d-3}$ if and
only if
$A$ is ample and $A$ satisfies property $N_{d-1}$ if and only if $A$ is ample
different from $-K_X$ ($-K_X$ satisfies $N_{d-3}$ but not $N_{d-2}$). 

\item{--} If $X=\bold F_e$, then  $A$ satisfies property $N_{e+1}$ if and
only if
$A$ is ample. 

\item{--} If $X=\bold P^2$, then  $A$ satisfies property $N_{0}$ if and
only if
$A$ is ample.

In addition, if $K_X^2=2$, then $A$ satisfies property $N_1$ if and only if $A$
is an ample line bundle different form $-K_X$.

\endproclaim

In the next result we show that very ampleness and projective normality are
equivalent for anticanonical surfaces:

\proclaim{\vapn } Let $X$ be an anticanonical surface. A line
bundle on $X$ is very ample if and only if it satisfies property
$N_0$. 
\endproclaim

\noindent {\it Proof.}
Let $L$ be very ample and let $C$ be smooth curve in $|L|$. Since $-K_X \otimes \Cal O_C$
is effective, $-K_X \cdot L \geq 3$, otherwise $L \otimes \Cal O_C$ would not be very
ample. Then $L$ satisfies property $N_0$ by \rsNp \ns. 
$\square$

\bigskip

We state and show now the result already announced dealing with line bundles
of the form $K_X+A_1+\cdots+A_n$, with $A_1,\dots,A_n$ ample. It follows
from \rsNp and Propositions 1.6, 1.9 and 1.10.

\proclaim{\Mukai } Let $X$ be an 
anticanonical surface. Let $A_1, \dots, A_n$ be  ample line
bundles on
$X$. If $n \geq p+4$, then $L=K_X +A_1 +\cdots +A_n$ satisfies property
$N_p$. This bound is achieved for $A_i=A$ and $(X,A)$ as in Example 1.15.

More precisely: 

\item{1)} If $X= \bold P^2$ and $n \geq \lceil \frac p 3 \rceil
+4$, then
$L$ satisfies property $N_p$. The bound is achieved for $A_i=\Cal O_{\bold
P^2}(1)$ and $p$ multiple of $3$. 

\item{2)} If $K_X^2=8$ and $n \geq \lceil \frac {p + 3} {4} \rceil + 2$, then
$L$ satisfies property $N_p$. The bound is achieved for $X = \bold F_0$, 
$A_i=C_0+f$ and $p \equiv 1 (4)$. More precisely, if
$X=\bold F_e$ and $n \geq \text{ max}(3, \lceil
\frac {p + 11} {e+4} \rceil)$, then $L$ satisfies property $N_p$. The bound is
achieved for instance for
$A_i=C_0+(e+1)f$ and for $p \equiv -11 (e+4)$. 

\item{3)} If $1 \leq K_X^2 \leq 7$ and $n \geq  \lceil \frac {p
+ 3} {K^2}
\rceil + 1$, then
$L$ satisfies property $N_p$. The bound is achieved for instance for $(X,A_i)$
as in Examples 1.13, 1.14 and 1.15, and for $p \equiv -3 (K_X^2)$.

\item{4)} If $1 \leq K_X^2 \leq 7$, $A_i \neq -K_X$, $A_i \neq -2K_X$ if
$K_X^2=1$ and
$n
\geq 
\lceil
\frac {p+K_X^2+3} {K_X^2+2} \rceil$, then $L$ satisfies property $N_p$. The
bound is achieved for $(X,A_i)$ as in Example 1.16, and for $p+K_X^2 \equiv -3
(K_X^2+2)$. 

\item{5)} If $2 \leq K_X^2 \leq 7$, $A_i\neq -K_X$, $A_i \neq -2K_X$ when
$K_X^2=2$, $A_i \neq -2K_X, -3K_X$ when $K_X^2=1$, it does not happen that 
$K_X+A_i$ is base-point-free, $A_i$ is very ample and $(X,A_i)$ is a conic
fibration under
$|K_X+A_i|$,  and
$n
\geq \text{ max}(2,
\lceil
\frac {p+K_X^2+3} {K_X^2+3} \rceil)$, then $L$ satisfies property $N_p$. The
bound is achieved for $(X,A_i)$ as in Example 1.17, and for $p+K_X^2 \equiv -3
(K_X^2+3)$.

\item{6)} If   $-1 \leq K_X^2 \leq 1$ and $n \geq p+3+K^2_X$, then
$L$ satisfies property $N_p$. The bound is sharp for $(X,A_i)$ as in  Examples
1.15, 1.18 and 1.19.

\item{7)} If   $K_X^2 \leq -2$ and $n \geq \text{ max}(2,p+3+K_X^2)$, then
$L$ satisfies property $N_p$. The bound is sharp for $(X,A_i)$ as in Examples
1.19  and 1.20.

\endproclaim

The next result we will prove is \ReiderNp \ns, which is an $N_p$ result in the 
same flavor of Reider's theorem for base-point-freeness and very ampleness. 
\ReiderNp shows that when $K_X^2 \geq 1$, a high self-intersection number
for $L$ implies by itself  property $N_p$ for $K_X + L$ for a large value of
$p$. This behavior is in contrast with the behavior that can be observed in
surfaces of Kodaira dimension $0$.

\proclaim{\ReiderNp }  Let $X$ be a rational  surface. Let
$L$ be a line bundle such that

\item{1)} $L \cdot C \geq 3$ for any  curve $C$ on $X$ and $L^2 \geq 10$ or

\item{1')} $K_X+L$ is very ample. 

\medskip

\noindent If $K_X^2 \geq 1$:

\item{2a)} Let $L^2 \geq (p+3)^2 +1$.

\medskip 

\noindent If $K_X^2 \geq 1$ and $L$ is not a multiple of $-K_X$ when $K_X^2=1$:

\item{2b)} Let $L^2 \geq (p+3)^2-1$. 

\medskip

\noindent If $K_X^2 \leq 0$:

\item{2c)} Let $-K_X \cdot L \geq p+3$.

\medskip
 
Then $K_X +L$ satisfies property $N_p$. 
\endproclaim

To prove \ReiderNp we will use
 \rsNp \ns, \vapn  and  the technical  lemma 1.25. The lemma is proven by
contradiction and the argument connects property  $N_p$ with termination of
adjunction. In particular it shows that a failure of certain property $N_p$ to
hold results in non termination of adjunction.

\proclaim{\term \ns} Let $X$ be a rational surface 
with
$K_X^2 >0$ and not isomorphic to $\bold P^2$. 
Let $L$ be a line
bundle on
$X$ such that 
  $K_X+L$ is
effective. Let $p \geq 1$ if $K_X^2 \leq 7$ and $p \geq 2$ if $K_X^2=8$. 
If
$L^2 \geq (p+3)^2 -1$ and $L$ is not a positive multiple of
$-K_X$ when $K_X^2=1$, then  
$-K_X
\cdot L
\geq p+3 + K_X^2$. 

\endproclaim

{\it Proof.}  
Assume that, $L$ is not a positive multiple of $-K_X$, or $K_X^2 \geq 2$. 
 We will prove the result by way of contradiction.
Assuming that $-K_X \cdot (K_X +L) \leq p+2$, we will prove by
induction that
$mK_X+L$ is effective for all $m \geq 2$. This contradicts the
termination of adjunction on a surface of Kodaira dimension
$-\infty$. 

If $m=1$, $K_X+L$ is effective by hypothesis. 
Let us see now that $2K_X+L$ is effective. Note that $L \neq -K_X$ because $L^2
\geq 15$ and $K_X^2 \leq 8$. Then, since $K_X + L$ is
effective,  $-K_X-L$ is not effective. Then by
Riemann-Roch
$$h^0(2K_X+L) \geq \frac 1 2 (2K_X+L)(K_X+L) + 1 \ .$$ Therefore we
need to see that $(2K_X +L)(K_X+L)
>  -2 $. This inequality is equivalent to 
$$L^2 +2K_X^2 +1 \geq 3(-K_X \cdot L) \ .$$
Since by assumption $-K_X\cdot L \leq p+2+K_X^2$, it suffices to check
that 
 $L^2 \geq 3(p+2)+K_X^2-1$. Then it
is enough to see that $(p+3)^2 \geq 3p+6+K_X^2$. This last
inequality is equivalent to $p^2+3p+3-K_X^2 \geq 0$, which holds for all
$p \geq 1$ if $K_X^2 \leq 7$ and for all $p \geq 2$ if $K_X^2=8$. 

Let $m \geq 2$. Now we assume $mK_X+L$ is effective, and we will show
$(m+1)K_X+L$ is effective. 
First consider the case when $K_X^2 =1$ and $L$ is not multiple of $-K_X$.  
Since $mK_X + L$ is
effective,  $h^0(-mK_X-L) =0$.   Then by Riemann--Roch
$$  h^0((m+1)K_X+L) \geq \frac 1 2 ((m+1)K_X +L)(mK_X+L) + 1  \ .$$ 
Thus we need to see that
$ ((m+1)K_X +L)(mK_X+L) \geq -1 $. 
 This is equivalent to 
$$L^2 +m(m+1)K_X^2 +1 \geq  (2m+1)(-K_X \cdot L) \ .$$
Since by assumption $-K_X\cdot L \leq p+2+K_X^2$, it suffices to check
that 
$$ L^2 \geq (2m+1)(p+2)+(-m^2+m+1)K_X^2 -1 \ . \eqno (1.25.1) $$
Recall that $K_X^2=1$. Then, in order to check 
(1.25.1) it is enough to show
that
$$(p+3)^2 -1\geq (2m+1)(p+2)-m^2+m \ .$$
This last inequality is equivalent to 
$$(p-m)^2 + 5(p-m) + 6 \geq 0 \ .$$ This holds for all integers $p$ and 
$m$.

Now we assume that $K_X^2 \geq 2$. Since $mK_X + L$ is
effective,  $h^0(-mK_X-L) \leq 1$.   Then by Riemann--Roch
$$  h^0((m+1)K_X+L) \geq \frac 1 2 ((m+1)K_X +L)(mK_X+L)   \ .$$ 
Thus we need to see that $ ((m+1)K_X +L)(mK_X+L) > 0 $. 
 This is equivalent to 
$$L^2 +m(m+1)K_X^2 > (2m+1)(-K_X \cdot L) \ .$$
Since by assumption $-K_X\cdot L \leq p+2+K_X^2$, it suffices to check
that 
$$ L^2 > (2m+1)(p+2)+(-m^2+m+1)K_X^2 \ . \eqno (1.25.2)$$
Since $m \geq 2$, $-m^2+m+1 < 0$. Then, in order to check 
(1.25.2) it is enough to show
that
$$(p+3)^2-1 > (2m+1)(p+2)-2(m^2-m-1) \ .$$
This last inequality is equivalent to 
$$p^2 + (5-2m)p + (2m^2-6m+4) > 0 \ .$$ This holds for any integer $p \geq 1$
and any
$m \geq 2$.

Summing up, we have just shown that,  under
the   hypothesis of the proposition and with the 
additional assumption that $-K_X
\cdot L
\leq p+2 + K_X^2$, 
$mK_X+L$ is effective for all  $m \geq 2$.  As pointed out before, 
this is a contradiction. Therefore
$-K_X
\cdot L
\geq p+3 +K_X^2$, as wished. $\square$

\bigskip
\noindent {\bf Remark 1.26.} \term  follows  from Hodge Index
Theorem  for some values of $p$ and $K_X^2$ but not for all.

\bigskip

\noindent (1.27) {\it Proof of \ReiderNp.}
By hypothesis 1') or by 1)  and Reider's Theorem it follows that
$K_X+L$ is very ample. Therefore by
\vapn   $K_X+L$ satisfies property $N_0$. 
 If $X = \bold P^2$, it follows clearly from 2b) that $-K_X \cdot (K_X+L) \geq
p+3$ if $p \geq 2$. Then by
\rsNp
$K_X+L$ satisfies property
$N_p$ (if $L$ is such that $L^2=16$, then one could say that $K_X+L$
satisfies $N_{\infty}$). We consider now $X=\bold F_e$. If $p \geq 2$, by \term
\ns, 2b) and \rsNp $L$ satisfies property $N_p$. If $p =1$ the result also
follows because if $K_X+L$ is very ample it can be seen that $L^2 \geq 18$.
   Now if $K_X^2=1$ and $L=m(-K_X)$, then
2a) implies $m \geq p+4$. Therefore it follows by \rsNp that $K_X+L$ satisfies
$N_p$. If $L$ is not a multiple of $-K_X$ or if $K_X^2 \geq 2$ and $X \neq
\bold P^2$, then it follows from 2b) and
\term that  $-K_X \cdot (K_X + L) \geq p+3$. Thus $K_X+L$ satisfies property
$N_p$ by \rsNp \ns.
 Finally if  
$K_X^2
\leq 0$,  by 2c) 
$-K_X
\cdot (K_X + L)  \geq p+3$. Thus $K_X+L$ satisfies 
property
$N_p$   by \rsNp \ns.
$\square$

\bigskip

\noindent{\bf Remark 1.28} \ReiderNp is optimal. To see the optimality when
2a) is assumed, take  
$X$ as in  Example 1.15 and $L=(p+4)(-K_X)$. To see the optimality when 1') and
2b) are assumed take $A$ as in Example 1.16, $e=0$, $l=6$ and $L=-K_X+A$. Then
$K_X+L$ satisfies property $N_1$ but not $N_2$ by \rsNp and $L^2=16$. On the
other hand if
$K_X^2
\leq 0$, assuming only hypothesis 1) and 2a) or 1') and 2a) do not suffice.
This can be seen taking
$(X,A)$ as in Example 1.18 and $L=3A$. Indeed $K_X+L$ satisfies property $N_0$
but not $N_1$; however, $L^2=27 > (2+3)^2 +1$. 

\bigskip

We end this section by showing the relation between the property $N_p$
satisfied by a line bundle $L$ and the termination of ampleness for
$mK_X+L$.

\proclaim{Theorem 1.29} Let $X$ be an anticanonical surface and let $L$ be a
line bundle satisfying property $N_p$ but not property $N_{p+1}$. 

\item{a)} If $X=\bold P^2$ and $m > \frac{p} {K_X^2}$;

\item{b)} if $X=\bold F_e$ and $m > \frac {p-e-1} {K_X^2}$;

\item{c)} if
$1 \leq K_X^2 \leq 7$, and
$m > \frac {p+3} {K_X^2}
-1$;

\item{d)} if
$1 \leq K_X^2 \leq 7$, $L$ is not a multiple of $-K_X$ and
$m > \frac {p+1} {K_X^2} -1$;

\item{e)} if $K_X^2 <0$, and $m < \frac {p+2} {K_X^2}$,

then $mK_X+L$ is not ample.
\endproclaim

\noindent {\it Proof.}
We outline the proof of c), d), e); a) and b) are similar. Let $1 \leq K_X^2
\leq 7$ and assume
$mK_X+L$ is ample. Then
$-K_X\cdot (mK_X+L)
\geq K_X^2$  by Proposition 1.10. This implies $-K_X\cdot L \geq (m+1)K_X^2$.
Now if
$L$ satisfies $N_p$ but not $N_{p+1}$,
$-K_X
\cdot L = p+3$. Then $(m+1)K_X^2 \leq p+3$, and $m \leq \frac {p+3} {K_X^2}
-1$. If $L$ is not a multiple of $-K_X$, then $-K_X\cdot (mK_X+L)
\geq K_X^2+2$ and we argue similarly. 

Let now $K_X^2 <0$ and assume again that $mK_X+L$ is ample.  Then $-K_X\cdot
(mK_X+L) \geq 1$. This implies $-K_X\cdot L \geq mK_X^2+1$. Again,  if
$L$ satisfies $N_p$ but not $N_{p+1}$,
$-K_X
\cdot L = p+3$. Then $mK_X^2 + 1 \leq p+3$, and $m \geq \frac {p+2} {K_X^2}$.
$\square$

\bigskip

\noindent{\bf Remark 1.30.} The bounds for $m$ in the previous theorem are
sharp. That is clear for the bound in a). For the bound in b) take
$L=C_0+(e+1)f$. For the bound in c) take $L=n(-K_X)$. For the bound in d) take
$L=A$, where $A$ is as in Example 1.16.

\heading 2. Fano $n$-folds of index greater than or equal to $n-1$
\endheading

In this section and in the next we prove results on syzygies of certain Fano
$n$-folds.   The first
attempt one would try to make to tackle this problem is to imitate the
arguments we carried on in Section 1. Given a Fano $n$-fold $X$ and a very ample
line bundle $L$, $H^1(rL)=0$ for all $r \geq 0$. Then  one could try to obtain
information about the free resolution of the image of $X$ from whatever
information is available on the free resolution of its general hyperplane
section $X'$, which is now of dimension
$n-1$. If no information is readily available on the resolution of $X'$, one
would iterate the argument, and taking successive hyperplane sections one could
read the Betti numbers of the resolution of $X$ from the resolution of a
surface or even, of a curve. This worked very well in the case of rational
surfaces because we ended with a curve $C$ and a line bundle $L_C$ on $C$ of
relatively high degree, and because of this high degree, we knew relevant
information about the resolution of the embedding of $C$ by $|L_C|$, thanks to
the results of Green and Lazarsfeld. However when the dimension is higher
we lose control of the hyperplane sections of $X$, or rather, there is much less
information available on the syzygies of the hyperplane sections. 

\medskip

We will
look at one example to explain what we mean. Let $X$ be $\bold P^n$ and
let $L=\Cal O_{\bold P^n}(d)$. If we consider the  intersection
of $n-2$ general divisors of $|L|$ we end with a surface in $\bold P^n$
which is a $(d,\dots,d)$ complete intersection. The only surfaces among those
which are rational surfaces (in fact, anticanonical) other than linear
$\bold P^2$ are the quadric and the cubic hypersurface in $\bold P^3$ and a
$(2,2)$ complete intersection in
$\bold P^4$. This means that, using essentially the same ideas as in Section
1, one is able to give a result as precise as
\rsNp regarding the property $N_p$ for $X=\bold P^3$ and $L=\Cal O_{\bold
P^2}(2), \Cal O_{\bold
P^2}(3)$ and for $X=\bold P^4$ and $L=\Cal O_{\bold
P^2}(2)$. Precisely one has that $\Cal O_{\bold P^3}(2)$ satisfies property
$N_5$ but not
$N_6$,
$\Cal O_{\bold P^3}(3)$ satisfies property $N_6$ but not $N_7$ and
$\Cal O_{\bold P^4}(2)$ satisfies property $N_5$ but not $N_6$. These
particular cases were known  (cf. \JPW and \OP \ns).  
However for other pairs $(X,L)$ the game turns out to be much more complicated.
For instance, if $(X,L)= (\bold P^3, \Cal O_{\bold P^3}(4))$ or $(\bold P^5,
\Cal O_{\bold P^5}(2))$ by the above process we arrive at a K3 surface. Then
knowing the syzygies of a K3 surfaces is equivalent to knowing the syzygies
of its hyperplane section, a canonical curve. On this much less information is
known and moreover, this information would depend (at least conjecturally) on
the Clifford index of the hyperplane section. Finally all the other complete
intersection surfaces of type $(d, \dots,d)$ are surfaces of general type.
Then $L_C$ will be a line bundle on the hyperplane section $C$ of degree
strictly less than $2g(C)-2$, and for those line bundles our knowledge is
even more incomplete than for the canonical line bundle.

\medskip
  
Because of all the above we need to carry out different arguments in this and in
the following section, but before that we will state   Theorem 2.1, which is
based upon the work done in Section 1.
Theorem 2.1 gives a necessary and sufficient condition for the line
bundle $H$ giving the index of the Fano $n$-fold of index $n-1$ to
satisfy property
$N_p$. Since a Fano $n$-fold $(X,H)$ of index $n+1$ or $n$ is $(\bold P^n,
\Cal O_{\bold P^n}(1))$ or a quadric hypersurface, we do not consider Fano
$n$-folds of these indices in the theorem:

\proclaim{\Fanopn} Let $X$ be a Fano $n$-fold. Assume there exists  
 an ample and base-point-free line bundle $H$ such that 
$K_X^*= (n-1)H$ (e.g., if $X$ is a Fano $n$-fold of index $n-1$). Assume
furthermore that
$H^n \geq p+3$. Then 
$H^1(M_{rH}\otimes lH)=0$ for all $r,l \geq 1$ and $H$ satisfies
property
$N_p$. 

\endproclaim

\noindent {\it Proof.} The result is proven by induction on the dimension $n$
of 
$X$, starting the induction in dimension $2$. If $n=2$, the results
follows from \rsNp \ns. Indeed, the only thing to be checked is
that $-K_X \cdot H \geq p+3$, $p \geq 0$. Then $-K_X=H$, so we have the required
inequality by hypothesis. 

Let us assume the result to be true for all dimensions from $2$ to
$n-1$, and we will prove it for $X$ of dimension $n$. We first see that $H$
satisfies property $N_0$. It suffices to prove that 
$$ H^0(rH) \otimes H^0(H) @>\alpha >> H^0((r+1)H) $$ surjects for all $r \geq
1$. Let $Y$ be a smooth irreducible member of $|H|$. By Kodaira
Vanishing Theorem,
$H^1(lH)=0$ for all
$l
\geq 0$, hence by \restobs \ns, it suffices to see that
  
$$0 \longrightarrow H^0(rH_Y) \otimes H^0(H_Y) @>\beta >>
H^0((r+1)H_Y)
\longrightarrow 0$$ surjects for all $r \geq 1$.

The variety $Y$ is a Fano $(n-1)$-fold, $K_Y=(K_X+H)_Y=-(n-2)H_Y$ and
$H_Y^{n-1}=H^n \geq 3$. Therefore
$\beta$ surjects by induction.

Therefore we have just proven that $H$ satisfies property $N_0$. To see it does
satisfies property $N_p$, we argue as in \rsNp \ns. Since $H^1(lH)=0$ for all
$l \geq 0$, $H$ satisfies the same property $N_p$ as $H_Y$. Then by
induction, since $H_Y^{n-1}=H^n \geq p+3$, $H$ satisfies property $N_p$. 
$\square$

\bigskip

\noindent{\bf Remark 2.2.} \Fanopn is in fact a characterization of the property
$N_p$ satisfied by $H$. This follows because the same is true for the
successive hyperplane section, which is a rational surface as we explained in
the proof of \Fanopn \ns. 

\bigskip

Now we want to show a result about the syzygies associated to the multiples of
the line bundle $H$. As pointed out at the beginning of this section we need to
use different ideas than those used in Section 1. 
Any result on the graded Betti numbers of the resolution of a variety can be
realized in terms of Koszul cohomology. This was shown by M. Green. For a
base-point-free line bundle $L$ we define the vector bundle $M_L$ as 
$$ 0 \longrightarrow M_L \longrightarrow H^0(L) \otimes \Cal O
\longrightarrow L \longrightarrow 0 \ . \eqno \seq $$ 
Then regarding property
$N_p$ Green showed the following criterion:

\proclaim {\GLlemma (\EL
\ns, Section 1.)} Let $L$ be an ample, globally generated line
bundle on a variety $X$. If the group 
\text{$H ^1(\bigwedge ^{p'+1} M_L \otimes
sL)$} vanishes for all $0 \leq p' \leq p$ and all $s \geq
1$, then
$L$ satisfies the
 property
$N_p$. If in addition  $H^1(rL) = 0$, for all $r \geq 1$,
then the above  is a necessary and sufficient condition
for $L$ to satisfy property $N_p$. 
\endproclaim

According to \GLlemma the results
(see for instance \rsNp) obtained in Section 1 for the syzygies of a
rational surface $X$ embedded by $L$ are equivalent to the vanishing of
$H^1(\wedge^{p'+1}M_L
\otimes sL)$ for all $0 \leq p' \leq p$ and $s \geq 1$. 
To carry out the  arguments for Fano
$n$-folds, we would need however to have the vanishing of
$H^1(M_L^{\otimes p+1}
\otimes sL)$, which does not follow in general from the vanishing of 
$H^1(\wedge^{p+1}M_L \otimes sL)$. That is the reason why we prove
\tensNp \ns. Before that, we need to state two auxiliary lemmas: 

\bigskip

\noindent {\bf \obsfree \ns.}
Let $E$ and $L_1, \dots , L_r$ be 
coherent sheaves on a variety 
$X.$
Consider the map $H^0(E) \otimes H^0(L_1
+ \cdots + L_r) @> \psi >> H^0(E \otimes
L_1 + \cdots + L_r)$ and the maps 
$$\displaylines{H^0(E) \otimes H^0(L_1
) @> \alpha_1 >> H^0(E \otimes L_1),
\cr
H^0(E\otimes L_1) \otimes H^0(L_2
) @> \alpha_2 >> H^0(E \otimes L_1 + L_2),\cr
 \dots ,
\cr
 H^0(E \otimes L_1 + \cdots + L_{r-1}) \otimes
H^0(L_r) @> \alpha_r >> H^0(E \otimes
L_1 + \cdots + L_r) \ .}$$
If $\alpha_1, \dots , \alpha_r$ are surjective
then $\psi$ is also surjective.

\proclaim{\Splemma (\GPfour \ns, Lemma 2.9)} Let $X$ be a projective variety,
let
$q$ be a nonnegative integer and let $F_i$ be a base-point-free line bundle
on $X$ for all $1 \leq i \leq q$. Let $Q$ be an effective line
bundle on $X$ and let $\frak q$ be a reduced and irreducible member
of
$|Q|$. 
Let $R$ be a line bundle and $G$ a sheaf on $X$ such that
\vskip .1 cm
\item{1.} $\text H^1(F_i \otimes Q^*)=0$
  \item {2.}  
$
\text H
^0(M_{(F_{i_1}
\otimes
\Cal O_\frak q)}
\otimes \dots \otimes M_{(F_{i_{q'}} \otimes \Cal O_\frak q) }\otimes
R \otimes \Cal O_\frak q)
\otimes \text H^0(G) \to$
\item{}
$\to \text H ^0(M_{(F_{i_1} \otimes 
\Cal O_\frak
q)}
\otimes \dots \otimes M_{(F_{i_{q'}}  \otimes \Cal O_\frak
q)}
\otimes R \otimes G \otimes \Cal O_\frak q) \ \text {surjects for all}
\  0 \leq q' \leq q$.

\noindent Then, for all $0 \leq q'' \leq q$ and any subset
$\{j_k\}
\subseteq \{i\}$ with $\#\{j_k\}=q''$ and for all $0 \leq k'
\leq q''$, 
$$\displaylines {\text H^0(M_{F_{j_1}} \otimes \dots \otimes
M_{F_{j_{k'}}}
\otimes M_{(F_{j_{k'+1}} \otimes \Cal O_\frak q)} \otimes \dots
\otimes M_{(F_{j_{q''}}\otimes \Cal O_\frak q)} \otimes R
\otimes \Cal O_\frak q) \otimes \text H^0(G) \to \cr
\text
H^0(M_{F_{j_1}} 
\otimes \dots \otimes M_{F_{j_{k'}}}
\otimes M_{(F_{j_{k'+1}} \otimes \Cal O_\frak q)} \otimes \dots
\otimes M_{(F_{j_{q''}} \otimes \Cal O_\frak q)} \otimes G \otimes R
\otimes \Cal O_\frak q)}$$ surjects.
\endproclaim

Now we are ready to prove

\proclaim{\tensNp} Let $X$ be a rational surface. Let $B$ be an ample and
base-point-free  line bundle such that
$-K_X
\cdot B \geq 4$ or $(X,B)=(\bold P^2,\Cal O_{\bold P^2}(1))$. Then
$H^1(M^{\otimes p+1}_{rB}\otimes lB)=0$ for all
$r
\geq 1$, $l \geq p$ and $p \geq 1$. In particular $lB$ satisfies
property $N_p$ for all $l \geq p$. 
\endproclaim

\noindent {\it Proof.}
We first do the case $-K_X \cdot B \geq 4$.    
The proof goes by induction on $p$. We start proving the case $p=1$.
We want to show that $H^1(M^{\otimes 2}_{rB}\otimes lB)=0$ for all $r
\geq 1$ and all $l \geq 1$. We tensor \seq associated to $rB$ by
$M_{rB} \otimes lB$ and take global sections. As a piece of the long
exact sequence of cohomology we obtain:

$$\displaylines{H^0(M_{rB} \otimes lB) \otimes H^0(rB) @>\alpha >> H^0(M_{rB}
\otimes (r+l)B) \cr 
\longrightarrow H^1(M^{\otimes 2}_{rB}\otimes lB)
\longrightarrow H^1(M_{rB} \otimes lB) \otimes H^0(rB) \ .}$$

By \genpn \ns, the last term of the above sequence is $0$. Then the
vanishing of $H^1(M^{\otimes 2}_{rB} \otimes lB)$ is equivalent to the
surjectivity of $\alpha$. To have the surjectivity of $\alpha$, by
\obsfree it suffices to see the surjectivity of

$$ H^0(M_{rB} \otimes lB) \otimes H^0(B) @>\beta >> H^0(M_{rB}
\otimes (l+1)B) \ .$$

Since $B$ is base-point-free and ample, we can choose $C$ smooth and
irreducible in $|B|$. Then to see the surjectivity of
$\beta$, by
\restobs and
\Spanishlemma
\ns, it suffices to see the surjectivity of

$$ H^0(M_{rB_C} \otimes lB_C) \otimes H^0(B_C) @>\gamma >>
H^0(M_{rB_C}
\otimes (l+1)B_C) \ .$$

To obtain the surjectivity of $\gamma$ we 
show
that the cokernel of $\gamma$ vanishes. The cokernel of $\gamma$ is
$H^1(M_{rB_C}^{\otimes 2} \otimes lB_C)$. Since $-K_X \cdot B \geq
4$, deg$B_C \geq 2g(C)+2$, then by
\Bu \ns, Theorem 1.12, $M_{rB_C}$ is semistable, and by \Mi \ns, Corollary 3.7
so is
$M_{rB_C}^{\otimes 2} \otimes lB_C$. A simple calculation shows that
$\mu(M_{rB_C}^{\otimes 2}
\otimes lB_C) > 2g(C)-2$ and we get the desired vanishing.

\medskip

Now we assume the result to be true for $1, \dots, p-1$ and we will
prove it for $p$, i.e., we want to see that $H^1(M^{\otimes
p+1}_{rB}\otimes lB)=0$. We go over the steps given to prove the case
$p=1$.  Since by induction $H^1(M^{\otimes p}_{rB}\otimes lB)=0$,
the vanishing we seek is equivalent to the surjectivity of 
 
$$H^0(M_{rB}^{\otimes p} \otimes lB) \otimes H^0(rB) @>\alpha >>
H^0(M_{rB}^{\otimes p}
\otimes (r+l)B) \ .$$

By \obsfree it suffices to see that
$$ H^0(M_{rB}^{\otimes p} \otimes lB) \otimes H^0(B) @>\beta >>
H^0(M_{rB}^{\otimes p}
\otimes (l+1)B) $$ surjects. Now choosing an irreducible and
smooth curve $C$ in $|B|$ and using
\restobs and
\Splemma
\ns, 

we see that it is enough to show that

$$ H^0(M_{rB_C}^{\otimes p} \otimes lB_C) \otimes H^0(B_C) @>\gamma
>> H^0(M_{rB_C}^{\otimes p}
\otimes (l+1)B_C)$$ surjects.

Finally $\gamma$ surjects if $H^1(M_{rB_C}^{\otimes p+1} \otimes
lB_C)=0$. By \Bu \ns, Theorem 1.12   and \Mi \ns, Corollary 3.7, the bundle
$M_{rB_C}^{\otimes p+1}
\otimes lB_C$ is semistable. On the other hand its slope $\mu$ is
$$\frac{-(p+1)rB^2}{rB^2-g(C)}+lB^2 \ ,$$ so we will conclude our
argument if we see that $\mu > 2g(C)-2$. This follows from $l \geq p
\geq 2, -K_X \cdot B \geq 4$. 

As in the case $p=1$, one can alternatively deduce the surjectivity of
$\gamma$ from \Bu \ns, Proposition 2.2.

To finish the proof we take care of the case $(X,B)=(\bold P^2, \Cal
O_{\bold P^2}(1)$. The proof goes along the same lines as before. We want the
vanishing of $H^1(M_{rB}^{\otimes p+1} \otimes lB)$ for all $r \geq
1$ and all $l \geq p$. This follows from the vanishing of $H^1(M_{\Cal
O_{\bold P^1(r)}}^{\otimes p+1} \otimes \Cal O_{\bold P^1}(l))$, which can be
easily checked as $M_{\Cal
O_{\bold P^1(r)}}$ is direct sum of copies of $\Cal O_{\bold P^1}(-1)$.
$\square$

\bigskip

We use  \tensNp to obtain the following results on Koszul cohomology and
syzygies for multiples of $H$.

\proclaim{\FanoNp} Let $X$ be a Fano $n$-fold. Assume there exists  
 an ample and base-point-free line bundle $H$ such that 
$-K_X= mH$, with $m \geq
n-1$ (e.g., if $X$ is a Fano $n$-fold of index $m \geq n-1$). Assume
furthermore that $H^n \geq 4$ if $m=n-1$. 
Then $H^1(M^{\otimes p+1}_{rH}\otimes lH)=0$ for
all
$r
\geq 1$, $l \geq p$ and $p \geq 1$.  
\endproclaim

\noindent {\it Proof.} 
To prove $H^1(M^{\otimes p+1}_{rH} \otimes lH)=0$  we argue by
induction on the dimension of
$X$, the cornerstone being now \tensNp \ns, and by induction on $p$.

Let first $p=1$.  We want to obtain the vanishing of
$H^1(M_{rH}^{\otimes 2} \otimes lH)$ and, as announced, we do it by
induction on the dimension $n$. If $n=2$ the result is a particular
 case of \text{\tensNp \ns.} In fact, $-K_X  = mH$ with $m \geq 1$. If
$m=1$, $-K_X \cdot H \geq 4$ follows directly by hypothesis. If $m \geq 2$,
$-K_X \cdot H \geq 2$, since $H$ is ample, and $-K_X \cdot H \geq 4$
unless $(X,H) = (\bold P^2, \Cal O_{\bold P^2}(1))$.

Now we assume the result to be
true for dimensions
$2$ to
$n-1$ and we will prove it for dimension
$n$. From
\seq we obtain
$$\displaylines{H^0(M_{rH} \otimes lH) \otimes H^0(rH) @>\alpha >> H^0(M_{rH}
\otimes (r+l)H) \cr 
\longrightarrow H^1(M_{rH}^{\otimes 2} \otimes lH)
\longrightarrow H^1(M_{rH} \otimes lH) \otimes H^0(rH) \ .}$$

By the vanishing of the first cohomology of the multiples of $H$,
\Fanopn implies the vanishing of the last term of the above sequence. 
Then the vanishing we seek is equivalent to the surjectivity of
$\alpha$. By \obsfree we see that it will suffice to prove that 
$$H^0(M_{rH} \otimes lH) \otimes H^0(H) @>\beta >> H^0(M_{rH} \otimes
(l+1)H)$$
surjects for all $r,l \geq 1$. The conditions needed to apply 
\restobs and
\Splemma 
follow from \Fanopn and from the fact that $H^1(M_{rH})=0$, for
$H^1(\Cal O_X)=0$. Then it suffices to see 
$$H^0(M_{rH_Y} \otimes lH_Y) \otimes H^0(H_Y) @>\gamma >> H^0(M_{rH_Y}
\otimes (l+1)H_Y)$$ surjects, which follows because the result is
true for $p=1$ and $Y$ of dimension $n-1$, by induction hypothesis. 

Now we assume the result to be true for $1,\dots, p-1$ and we will
prove it for $p$. We argue again by induction on the dimension. In
dimension $2$ the result  follows from \tensNp as in the case $p=1$.
We therefore assume the result to be
true for $p$ and dimensions
$2$ to
$n-1$ and we will prove it for $p$ and dimension
$n$, i.e., we will prove $H^1(M^{\otimes p+1}_{rH}\otimes lH)=0$ for
all
$r
\geq 1$, $l \geq p$.
From
\seq we obtain
$$\displaylines{H^0(M_{rH}^{\otimes p} \otimes lH) \otimes H^0(rH) @>\alpha >>
H^0(M_{rH}^{\otimes p}
\otimes (r+l)H) \cr 
\longrightarrow H^1(M_{rH}^{\otimes p+1} \otimes lH)
\longrightarrow H^1(M_{rH}^{\otimes p} \otimes lH) \otimes H^0(rH) \
.}$$

By induction hypothesis on $p$,
we have the vanishing of the last term of the above sequence. 
Therefore the vanishing we seek is equivalent to the surjectivity of
$\alpha$. By \obsfree we see that it will suffice to prove that 
$$H^0(M_{rH}^{\otimes p} \otimes lH) \otimes H^0(H) @>\beta >>
H^0(M_{rH}^{\otimes p} \otimes (l+1)H)$$
surjects for all $r \geq 1,l \geq p$. The conditions needed to apply 
\restobs and
\Splemma 

follow from induction hypothesis on $p$. Then it suffices to see 
$$H^0(M_{rH_Y}^{\otimes p} \otimes lH_Y) \otimes H^0(H_Y) @>\gamma >>
H^0(M_{rH_Y}^{\otimes p}
\otimes (l+1)H_Y)$$ surjects, which follows because the result is
true for $p$ and $Y$ of dimension $n-1$, by induction hypothesis on
$n$.
$\square$

\proclaim{Corollary 2.8} Let $X$ be a Fano $n$-fold. Assume there 
exists  
 an ample and base-point-free line bundle $H$ such that 
$-K_X= mH$, with $m \geq
n-1$ (e.g., if $X$ is a Fano $n$-fold of index $m \geq 1$). Assume
furthermore that $H^n \geq 4$ if $m=n-1$. Then $lH$ satisfies
property $N_p$ for all $l \geq p$. 
\endproclaim

\noindent {\it Proof.} 

Since we work in characteristic $0$, then 
$$H^1(\bigwedge^{\otimes i} M_{lH}
\otimes slH)=0, \text{ for all }
l \geq p, s \geq 1, 1 \leq i \leq p+1 \ .$$  
Therefore, according to
\GL \ns,
$lH$ satisfies property
$N_p$.

\bigskip

\heading 3. Fano $n$-folds of index $n-3$. \endheading

In this last section we deal with Fano $n$-folds of index $n-3$, $n \geq
4$. We use the same ideas as in the second part of Section 2. We start by
studying under what conditions $mH$ is very ample and satisfies property $N_0$:

\proclaim{\CYpn}
Let $X$ be a Fano $n$-fold such that $-K_X=mH$,
$H$ is ample and base-point-free, and $m=n-3 \geq 1$ (for instance, if 
$X$ is a Fano $n$-fold of index $m=n-3$). Let
$L=kH$.    

\item{(1)} If $k \geq 4$, then $L$ satisfies property $N_0$. 
\item{(2)} Let $k=3$; $L$ satisfies property $N_0$  if and only if the morphism
induced by $|H|$ does not map $X \ 2:1$ onto $\bold P^n$. 
\item{(3)} Let $k=2$; if
 $|H|$ does not map $X$ onto a variety of minimal
degree other than $\bold P^n$ nor  maps $X \ 2:1$ onto 
$\bold P^n$, then $L$ satisfies property $N_0$.   
\endproclaim

{\it Proof.}
We will prove the result by induction on $n$. We start at $n=4$. We
have to deal with several cases. 
\medskip

{\it Case 1: $k \geq 4$.}
We want to prove that
$L=kH$ satisfies property $N_0$, or equivalently,  that 
$$H^0(skH) \otimes H^0(kH) \longrightarrow H^0((s+1)kH) $$ surjects for
all $s \geq 1$. 
This follows from a  more general result, namely,
$$H^0(kH) \otimes H^0(H) @>\alpha >> H^0((k+1)H) $$
surjects for all $k \geq 4$. Indeed, the surjectivity of $\alpha$  follows
from
\CM \ns,  p. 41, Theorem 2, since by Kodaira Vanishing Theorem, $H$ is
$3$-regular. 
\medskip

{\it Case 2: $k=3$.} We distinguish  two
subcases. 

{\it Case 2.1:} First let us assume  that $|H|$ does not map $X$ onto
$\bold P^4$. We will show $$H^0(kH) \otimes H^0(H) @>\alpha >>
H^0((k+1)H) $$ for all $k \geq 3$.  If $k \geq 4$, we have already seen
that $\alpha$ surjects. To show the surjectivity for $k=3$, let
$Y$ be a smooth irreducible member of
$|H|$. By
\restobs
\ns, since
$H^1(lH)=0$, 
 for all $l \geq 0$ by Kodaira Vanishing Theorem, it will suffice to
check that 
$$H^0(3H_Y) \otimes H^0(H_Y) @> \beta >> H^0(4H_Y) $$
surjects. By adjunction, $(Y,H_Y)$
is a polarized Calabi-Yau threefold, and, since $H^1(\Cal O_X)=0$, $|H_Y|$
does not map
$Y$ onto $\bold P^3$. Then $\beta$ surjects (cf. \GPthree \ns, proof of
Theorem 1.4, case 1).

\medskip
 
{\it Case 2.2:} Now assume  $|H|$ does  map $X$
onto $\bold P^4$. We assume first that the map induced by $|H|$ is not
$2:1$ and we will prove that 
$$H^0(3lH) \otimes H^0(3H) \longrightarrow H^0(3(l+1)H) $$ surjects for
all $l \geq 1$. For that, by \freeobs it is enough that 
$$\displaylines{H^0(rH) \otimes H^0(H) @>\alpha >> H^0((r+1)H) \text{
for all }r \geq 5 \cr
H^0(3H) \otimes H^0(2H) @>\gamma >> H^0(5H)}$$
surject. The map $\alpha$ was seen to be surjective in Case 1.
For the surjectivity of $\gamma$ we choose  smooth irreducible $Y$ in
$|H|$ and we write the following commutative diagram, 

$$\matrix
H^0(2H) \otimes  H^0(2H) & 
\hookrightarrow & H^0(2H) \otimes 
H^0(3H)&
\twoheadrightarrow & H^0(2H) \otimes 
H^0(3H_Y) \cr 
@VV\delta V @VV\gamma V 
@VV\epsilon V \cr
 H^0(4H) & \hookrightarrow &
H^0(5H)& \twoheadrightarrow &
H^0(5H_Y)  \ ,
\endmatrix
$$
obtained from the sequence
$$ 0 \longrightarrow H^* \longrightarrow \Cal O_X 
\longrightarrow \Cal O_Y
\longrightarrow 0
\eqno
(3.1.1)
\ , $$ having in account that $H^1(rH) = 0$ for all $r \geq 0$.  To see
the surjectivity of
$\delta$ we construct another diagram arising from $(3.1.1)$:

$$\matrix
H^0(2H) \otimes  H^0(H) & 
\hookrightarrow & H^0(2H) \otimes 
H^0(2H)&
\twoheadrightarrow & H^0(2H) \otimes 
H^0(2H_Y) \cr @VV\eta V @VV\delta V 
@VV\nu V\cr
 H^0(3H) & \hookrightarrow &
H^0(4H)& \twoheadrightarrow &
H^0(4H_Y)  \ .
\endmatrix
$$

Thus we need to see that $\eta$, $\epsilon$ and $\nu$ are surjective. By
\restobs \ns, for the surjectivity of $\eta$ it suffices to see the
surjectivity of 
$$H^0(2H_Y) \otimes H^0(H_Y) @> \mu >> H^0(3H_Y)$$
and for the surjectivity of $\epsilon$ and $\nu$, since
$H^1(H)=0$, it suffices to see that

$$\displaylines{ H^0(3H_Y) \otimes H^0(2H_Y) @> \pi >> H^0(5H_Y)\cr
H^0(2H_Y) \otimes H^0(2H_Y) @> \rho >> H^0(4H_Y)}$$
both surject. $(Y,H_Y)$ is a polarized Calabi-Yau threefold such that the
morphism induced by 
$|H_Y|$ maps $Y$ onto to $\bold P^3$ but is not $2:1$. Then the 
surjectivity of
$\pi$ and
$\rho$ is explicitly claimed (and proved) in the proof of Theorem 1.4,
\GPthree and the surjectivity of
$\mu$ is also proved there, although not explicitly said. 

To end Case 2.2 assume $|H|$ maps $X$ $2:1$ onto $\bold P^4$. If $Y$ is a
smooth and irreducible member of $|H|$, then $|H_Y|$ maps $Y$ $2:1$ onto
$\bold P^3$, so according to \GPthree \ns, Theorem 1.4, $H_Y$ is not very ample,
nor is $H$.
\medskip

{\it Case 3: $k=2$}. We want to see that if $H$ is ample, 
base-point-free and $|H|$ does not map $X$ onto a variety of minimal
degree (different from $\bold P^4$), nor does it map $X$ $2:1$ onto
$\bold P^4$, then 
$$H^0(2lH) \otimes H^0(2H) \longrightarrow H^0((2l+2)H)$$
surjects for all $l \geq 1$.

{\it Case 3.1:} Assume $h^0(H) \geq 6$.  
By \freeobs \ns, it suffices to prove that
$$H^0(rH) \otimes H^0(H) @> \alpha >> H^0((r+1)H)$$ surjects for all
$r \geq 2$. If $r \geq 3$, we have seen that $\alpha$ surjects while
proving Case 1 and Case 2.1. To see that
$$H^0(2H) \otimes H^0(H) @> \alpha >> H^0(3H)$$ surjects let $Y$ be as
before a smooth and irreducible member of $|H|$. By \restobs it suffices
to see that 
$$H^0(2H_Y) \otimes H^0(H_Y) @> \beta >> H^0(3H_Y) $$
surjects. Now $(Y,H_Y)$ is a polarized Calabi-Yau satisfying the
hypothesis of Theorem 1.7.1, \GPthree \ns. Then the map $\beta$ surjects, as
seen in the proof of Theorem 1.7, \GPthree \ns. 

{\it Case 3.2:} Assume now that $h^0(H) = 5$ and that the map induced by
$|H|$ from $X$ onto $\bold P^4$ is not $2:1$. By \freeobs it suffices to
see that 
$$\displaylines{H^0(rH) \otimes H^0(H) \longrightarrow H^0((r+1)H)
\text{ for all }r \geq 4 \text{ and }\cr  H^0(2H) \otimes H^0(2H)
\longrightarrow H^0(4H)}$$
surject, and this has been proved in Case 1 and Case 2.2. 

\medskip

To prove the result for a Fano variety of arbitrary dimension $n$ we
will argue by induction. Let $Y$ be smooth and irreducible member of
$|H|$. If $-K_X = (n-3)H$, then $-K_Y = (n-2)H_Y$ and
$h^0(H_Y)=h^0(H)-1$. Moreover, if the image of $X$ by the morphism
induced by $|H|$ is a variety of minimal degree, so is the image of
$Y$ by the morphism induced by $|H_Y|$, and the degree of both
morphisms is the same. Going over the arguments in the case
$n=4$, we see that the key point was to show the surjectivity of the
maps

$$\displaylines{  H^0(rH) \otimes H^0(H) @> \alpha >> H^0((r+1)H)
\text{ for all } k \geq 4 ,  \cr 
 H^0(3H) \otimes H^0(H) @> \beta >>
H^0(4H) \text{ if } h^0(H) \geq n+2
  \cr 
 H^0(2H) \otimes H^0(H) @>\gamma >> H^0(3H) \text{ if } h^0(H) \geq
n+2 \text{ and the image of $X$ } \cr 
\text{ by the map induced by $|H|$ is not 
a variety of minimal degree }\cr 
 H^0(3H)
\otimes H^0(2H) @> \delta >> H^0(5H) \text{ if } h^0(H) =
n+1 \text{ and the  map induced by $|H|$ is not }
2:1   \cr
 H^0(2H) \otimes H^0(2H) @>\epsilon >> H^0(4H) \text{ if } h^0(H) =
n+1 \text{ and the  map induced by $|H|$ is not }
2:1    \cr
 H^0(2H) \otimes H^0(H) @> \eta >> H^0(3H) \text{ if } h^0(H) =
n+1 \text{ and the  map induced by $|H|$ is not }
2:1  \ . } $$

Let us assume therefore the surjectivity of all the above maps
 for Fano varieties $(X,H)$ of dimension
$4,
\dots, n-1$. Arguing as in the case $n=4$, by \restobs and because
$H^1(rH)=0$ for all $r \geq 0$, we have that all the above maps are
also surjective if the dimension of $X$ is $n$.

Then by \freeobs the surjectivity of $\alpha$ implies (1), the surjectivity of
$\alpha$, $\beta$ and $\delta$ implies the ``if" part of (2), and the
surjectivity of
$\alpha$, $\beta$, $\gamma$ and $\epsilon$ implies (3). Finally, if the morphism
induced by $|H|$ is a double cover of $\bold P^n$, so is the morphism induced on
$Y$ by $|H_Y|$, when $Y$ is an irreducible and smooth member of $|H|$.
Then $3H_Y$ is not very ample (by induction hypothesis), nor is $3H$.   
\ $\square$

\bigskip

\proclaim{\CYNp}
Let $X$ be a Fano $n$-fold of index $m=n-3$. Assume that $-K_X=mH$,
and $H$ is ample and base-point-free. Let $L=kH$.
Assume furthermore that $h^0(H) \geq n+2$, i.e.,  that $|H|$ do
not map
$X$ onto
$\bold P^n$.  If
$k
\geq p+2$ and $p \geq 1$, then $L$ satisfies
property
$N_p$. 
\endproclaim

\noindent {\it Proof.} 
The proof is again by induction on the dimension. The first step is
$n=4$. Given $p \geq 1$, we want to prove that 
$$  H^1(M^{\otimes p+1}_L
\otimes sL)=0 \text{ if $L=kH$, $k \geq p+2$, $l \geq 1$. } \eqno (3.2.1) $$ 
 We will prove this more general
fact,  namely that 
$$  H^1(M^{\otimes p+1}_{kH}
\otimes sH)=0 \text{ for all  $k \geq p+2$, $s \geq p+2$.} \eqno (3.2.2) $$
We argue now by induction on $p$. If $p=1$, we want to prove that 
$$  H^1(M^{\otimes 2}_{kH}
\otimes sH)=0 \text{ for all $k \geq 3$, $s \geq 3$.} \eqno (3.2.3) $$
Since $H^1(M_{kH} \otimes sH) = 0$ by \CYpn \ns, then (3.2.3) is equivalent to
the surjectivity of
$$H^0(M_{kH}\otimes sH) \otimes H^0(kH) \longrightarrow H^0(M_{kH}\otimes
(k+s)H) \text{ for all $k \geq 3$, $s \geq 3$,}$$ and by \freeobs it
suffices to see that

$$H^0(M_{kH}\otimes H) \otimes H^0(H) \longrightarrow H^0(M_{kH}\otimes
(k+1)H) \text{ surjects for all $k \geq 3$.}$$

Choose smooth and irreducible $3$-fold $Y$ in $|H|$. Because of \CYpn
\ns,
\restobs and
\Splemma it suffices to have the surjectivity of

$$H^0(M_{kH_Y}\otimes H_Y) \otimes H^0(H_Y) \longrightarrow
H^0(M_{kH_Y}\otimes (k+1)H_Y) $$

for all $k \geq 3$. By adjunction $Y$ is a Calabi-Yau threefold and since
$H^1(\Cal O_X)=0$, $h^0(H_Y) \geq 5$, hence from the proof of
\GPthree \ns, Theorem 1.4, the above map is surjective.

We complete now the proof of the result for $n =4$. We may assume the
result proved until $p-1$. Now we want to see that 
$$  H^1(M^{\otimes p+1}_{kH}
\otimes sH)=0 \text{ for all  $k \geq p+2$, $s \geq p+2$. } \eqno (3.2.4) $$ 
We argue similarly to the case $p=1$. By induction on $p$ we may conclude
that the sought vanishing is equivalent to the surjectivity of

$$H^0(M_{kH}^{\otimes p}\otimes sH) \otimes H^0(kH) \longrightarrow
H^0(M_{kH}^{\otimes p} \otimes (k+s)H) \text{ for all $k \geq p+2$, $s \geq
p+2$,}$$ and by
\freeobs it suffices to see that

$$H^0(M_{kH}^{\otimes p}\otimes H) \otimes H^0(H) \longrightarrow
H^0(M_{kH}^{\otimes p}\otimes (k+1)H) \text{ surjects for all $k \geq p+2$.}$$

Finally we choose a smooth and irreducible $3$-fold $Y$ in $|H|$. Because
of
\CYpn
\ns,
\restobs and
\Splemma it suffices to have the surjectivity of

$$H^0(M_{kH_Y}\otimes H_Y) \otimes H^0(H_Y) \longrightarrow
H^0(M_{kH_Y}\otimes (k+1)H_Y) $$

for all $k \geq p+2$. By adjunction $Y$ is a Calabi-Yau threefold and
since
$H^1(\Cal O_X)=0$, $h^0(H_Y) \geq 5$, hence according to the proof of 
\GPthree \ns, Theorem 1.4, the above map is surjective.

Now assume $n > 4$. Recall that we want to show that  
$$  H^1(M^{\otimes p+1}_L
\otimes sL)=0 \text{ if $L=kH$, $k \geq p+2$, $l \geq 1$, } \eqno (3.2.5) $$ 
and as before we  will prove this more general
fact,  namely that 
$$  H^1(M^{\otimes p+1}_{kH}
\otimes sH)=0 \text{ for all  $k \geq p+2$, $s \geq p+2$.} \eqno (3.2.6) $$
We argue now by induction on $p$ and $n$. If $p=1$, we want to prove that 
$$  H^1(M^{\otimes 2}_{kH}
\otimes sH)=0 \text{ for all $k \geq 3$, $s \geq 3$.} \eqno (3.2.7) $$
Since $H^1(M_{kH} \otimes sH) = 0$ by \CYpn \ns, then (3.2.7) is equivalent to
the surjectivity of
$$H^0(M_{kH}\otimes sH) \otimes H^0(kH) \longrightarrow H^0(M_{kH}\otimes
(k+s)H) \text{ for all $k \geq 3$, $s \geq 3$,}$$ and by \freeobs it
suffices to see that

$$H^0(M_{kH}\otimes H) \otimes H^0(H) @>\alpha>> H^0(M_{kH}\otimes
(k+1)H) \text{ surjects for all $k \geq 3$.}$$ 

The surjectivity of this
map has been proven under the hypothesis of the theorem, when the
dimension
$n$ of
$X$ is
$4$, and we will assume it proved also if dimension of $X$ is $n-1$.

Choose smooth and irreducible $(n-1)$-fold $Y$ in $|H|$. Because of \CYpn
\ns,
\restobs and
\Splemma it suffices to have the surjectivity of

$$H^0(M_{kH_Y}\otimes H_Y) \otimes H^0(H_Y) \longrightarrow
H^0(M_{kH_Y}\otimes (k+1)H_Y) $$

for all $k \geq 3$. By adjunction $-K_Y=(n-4)H_Y$  and since
$H^1(\Cal O_X)=0$, $h^0(H_Y) \geq n+1$, hence by induction hypothesis,
$\alpha$ surjects.

The proof of the general case follows the same steps as the proof for
$n=4$. Recall that we want to proof that 
$$   H^1(M^{\otimes p+1}_{kH}
\otimes sH)=0 \text{ for all  $k \geq p+2$, $s \geq p+2$.} \eqno (3.2.8) $$
Using now the induction hypothesis for $n-1$ 
 one conclude
the result, exactly in the same fashion as we have just done when $p=1$.
$\square$

\heading References \endheading

\roster

\item"\BS" M.C. Beltrametti and A.J. Sommese, {\it The Adjunction Theory of
Complex Projective Varieties}, Walter de Gruyter, 1995.

\item"\Bu"
D. Butler, {\it Normal generation of vector bundles over
a curve}, J. Differential 
Geometry {\bf 39} (1994) 1-34.

\item"\EL" L. Ein and R. Lazarsfeld, {\it Koszul cohomology and syzygies of
projective
 varieties}, Inv. Math. {\bf  111} (1993), 51-67.

\item "\GPone" F.J. Gallego and B.P. Purnaprajna {\it Normal
presentation on elliptic ruled surfaces}, J.  Algebra {\bf 186}, 
(1996),
597-625.

\item"\GPtwo" \hbox{\leaders \hrule  \hskip .6 cm}\hskip .05 cm , {\it Higher
syzygies of elliptic ruled surfaces}, J.  Algebra {\bf 186}, (1996),
626-659.

\item"\GPthree" \hbox{\leaders \hrule  \hskip .6 cm}\hskip .05 cm , 
{\it Very ampleness and higher syzygies for Calabi-Yau threefolds},
Math. Ann. {\bf 312} (1998), 133-149.

\item"\GPfour" \hbox{\leaders \hrule  \hskip .6 cm}\hskip .05 cm , {\it
Projective normality and syzygies of algebraic surfaces}, J. reine
angew. Math. {\bf 506} (1999), 145-180.

\item"\GPfive" \hbox{\leaders \hrule  \hskip .6 cm}\hskip .05 cm , 
{\it Vanishing theorems and syzygies for K3 surfaces and Fano
varieties}, to appear in J. Pure   App. Alg.

\item"\Green" M. Green,
{\it Koszul cohomology and the geometry of projective varieties},
J. Differential Geometry {\bf
19} (1984), 125-171.

\item"\GL" M. Green and R. Lazarsfeld, 
 {\it Some results  on
the syzygies of finite sets and  
algebraic curves}, Compositio
Math. {\bf 67} (1989), 301-314. 

\item"\Hbone" B. Harbourne, {\it Anticanonical rational surfaces}, Trans. A. M.
S.  {\bf 349} (1997), 1191-1208.

\item"\Hbtwo" \hbox{\leaders \hrule  \hskip .6 cm}\hskip .05 cm ,
{\it Birational morphisms of rational surfaces}, J. Algebra {\bf 190} (1997),
145-162.

\item"\Ha" R. Hartshorne, Algebraic Geometry,
Springer--Verlag 1977.

\item"\Hoone" Y. Homma, {\it Projective normality and the
defining equations of ample 
invertible sheaves on elliptic ruled
surfaces with $e \geq 0$}, Natural Science 
Report, Ochanomizu
Univ. {\bf 31} (1980) 61-73.

\item"\Hotwo" 
\hbox{\leaders \hrule  \hskip .6 cm}\hskip .05 cm , {\it Projective
normality and the defining equations of an elliptic ruled 
surface with negative
invariant}, Natural Science Report, Ochanomizu
Univ. {\bf 33} (1982) 17-26.

\item"\JPW" T. Josefiak, P. Pragacz and J. Weyman, {\it Resolutions of
determinantal varieties and tensor complexes associated with symmetric and
antisymmetric matrix}, Asterisque 87-88 (1981), 109-189.

\item"\Mi" Y. Miyaoka, {\it The Chern class and Kodaira
dimension of a minimal variety}, 
Algebraic Geometry --Sendai 1985,
Advanced Studies in Pure Math., 
Vol. 10, 
North-Holland,
Amsterdam, 449-476.

\item "\CM" D. Mumford
{\it Varieties defined by quadratic equations}, Corso CIME in 
Questions on Algebraic Varieties,
Rome, 1970, 30-100.

\item"\OP" G. Ottaviani and R. Paoletti, {\it Syzygies of Veronese embeddings},
preprint.

\endroster

\enddocument

\end